\documentclass[11pt]{amsart}
\usepackage[dvips]{graphicx}
\usepackage{amsmath}
\usepackage{amsfonts}
\usepackage{amssymb}
\usepackage{psfrag, color}
\usepackage{subfigure}
\usepackage{lscape}
\usepackage{fancyhdr}
\usepackage{a4,epsfig}

\def\oh{\frac{1}{2}}
\def\fn{f_\nu}
\def\gn{\gamma_\nu}
\def\fns{f_\nu^\ast}

\def\fnpu{f_\nu(u^+_h(x,t))}

\def\Ojn{\Omega_j^n}
\def\OT{\Omega}
\def\Gin{\Gamma_{in}}
\def\Gout{\Gamma_{out}}
\def\Gjn{\Gamma_j^n}

\def\ab{\overline{a}(u;u_h)}
\def\ebkn{\bar\eta_k^n}

\def\vps{\varphi_\sharp}
\def\tol{T\!ol}
\def\cfl{C\!F\!L}
\newtheorem{theorem}{Theorem}
\newtheorem{remark}[theorem]{Remark}

\newtheorem{definition}[theorem]{Definition}

\begin{document}
\title{On adaptive timestepping for weakly instationary solutions of hyperbolic conservation laws via adjoint error control}
\author{Christina Steiner}
\author{Sebastian Noelle}

\begin{abstract}
We study a recent timestep adaptation technique for hyperbolic conservation laws.
The key tool is a space-time splitting of adjoint error representations for target functionals due
to S\"uli\cite{Sueli} and Hartmann\cite{Hartmann1998}. It provides an efficient choice of timesteps
for implicit computations of weakly instationary flows. The timestep will be very large in regions
of stationary flow, and become small when a perturbation enters the flow field.  Besides using
adjoint techniques which are already well-established, we also add a new ingredient which
simplifies the computation of the dual problem. Due to Galerkin orthogonality, the dual solution
$\varphi$ does not enter the error representation as such. Instead, the relevant term is the difference
of the dual solution and its projection to the finite element space, $\varphi-\varphi_h$. We can show that it
is therefore sufficient to compute the spatial gradient of the dual solution, $w=\nabla
\varphi$. This
gradient satisfies a conservation law instead of a transport equation, and it can therefore be
computed with the same algorithm as the forward problem, and in the same finite element space. We
demonstrate the capabilities of the approach  for a weakly instationary test problem for scalar
conservation laws.
\end{abstract}

\maketitle
\tableofcontents
\section{Introduction}

For explicit calculations of instationary solutions to hyperbolic
conservation laws, the timestep is dictated by the CFL condition due
to Courant, Friedrichs and Lewy \cite{CourantFriedrichsLewy1928},
which requires that the numerical speed of propagation should be at
least as large as the physical one. For implicit schemes, the CFL
condition does not provide a restriction, since the numerical speed of
propagation is infinite.  Depending on the equations and the scheme,
restrictions may come in via the stiffness of the resulting
nonlinear problem. These restrictions are usually not as strict as in
the explicit case, where the CFL number should be below unity. For
implicit calculations, CFL numbers of 10, 100 or even 1000 may well be
possible. Therefore, it is a serious question how large the timestep,
i.e. the CFL number, should be chosen. 

We are particularly interested in timestep control which is based upon 
computable, a-posteriori error estimates. In \cite{KroenerOhlberger2000,
Ohlberger2001} Kr\"oner and Ohlberger based their space-time adaptivity  upon
$L^1$, Kuznetsov type estimates for scalar conservation laws.  In
\cite{Johnson2,Johnson3,Johnson4,Johnson5,Johnson6}, Eriksson and Johnson developed space-time
adaptive methods for parabolic pde's. These a-posteriori error estimates require
the solution of an adjoint problem. A space-time projection of the adjoint
solution makes it possible to consider spatial and temporal error separately.
They closed the error estimates by an a-priori bound on the dual solution. In
\cite{Sueli, SueliHouston2003}, S\"uli and Houston developed an analogous
approach for hyperbolic transport equations.

The work  of Eriksson and Johnson has been extended by many authors, see, for
example, the review articles of Becker and Rannacher \cite{BeckerRannacher1996,
BeckerRannacher2000}
and of Hoffman and Johnson \cite{HoffmanJohnson}. We would like to mention that we learned
a lot about these developments from the unpublished thesis of Ralf Hartmann
\cite{Hartmann1998}. Instead of relying upon an (usually pessimistic) a-priori
error estimate for the adjoint solution, Hartmann and others
\cite{HartmannHoustonSISC2002, SueliHouston2003} {\em computed} the adjoint solution and hence
obtained an (in principle exact) error representation.

More recently these methods have also been developed
for hyperbolic problem by Barth, Hartmann, Houston, Giles, S\"uli,
Schwab and others. An excellent collection of review papers may
be found in \cite{BarthDeconinck2003}.

Let us briefly summarize the space-time splitting of the adjoint error
representation (see \cite{Johnson2,Johnson3,Johnson4,Johnson5,Johnson6,
BeckerRannacher1996,Sueli,Hartmann1998} for details). The error representation
expresses the error in a target  functional as a scalar product of the finite
element residual with the dual solution. This error representation is decomposed
into separate spatial and temporal components. The spatial part will decrease
under refinement of the spatial grid, and the temporal part under refinement of
the timestep. Technically, this decomposition is achieved by inserting an
additional projection. Usually, in the error representation, one subtracts from
the dual solution its projection onto space-time polynomials. Now, we also
insert the projection of the dual solution onto polynomials in time having
values which are $H^1$ functions with respect to space.

This splitting can be used to develop a strategy for a local choice of
timestep. Here we add to the
results in \cite{Sueli,Hartmann1998} by studying a weakly instationary
solution to Burgers' equation, for which the timestep will be very large 
(and we will quantify this) in regions of stationary flow, and become 
small when a perturbation enters the flow field. We believe that this type 
of flow is a prime example where the space-time splitting can become useful.

Besides applying adjoint techniques which are already well-established to a new
test problem, we also add a new ingredient which simplifies and accelerates the
computation of the dual problem. Due to Galerkin orthogonality, the dual
solution $\varphi$ does not enter the error representation as such. Instead, the
relevant term is the difference of the dual solution and its projection to the
finite element space, $\varphi-\varphi_h$. We can show that it is therefore sufficient to
compute the spatial gradient of the dual solution, $w=\nabla \varphi$. This gradient
satisfies a conservation law instead of a transport equation, and it can
therefore be computed with the same algorithm as the forward problem, and in the
same finite element space.

Our goal here is time step adaptation. Ultimately, this will become a
building block of an aerodynamic and aeroelastic solver which is
currently being developed by the SFB 401 research group at RWTH Aachen
\cite{BramkampLambyMueller2004}. In that solver, multiscale analysis
is used to compress data, coarsen and refine the spatial grid. Time stepping for
instationary problems is done by a methods of lines approach, using
explicit or implicit Runge-Kutta schemes. The latter is, of course, a
standard set-up used for aerodynamic, or conservation law, solvers.

In the aerodynamical applications which we have in mind, we may have
to resolve many different features of the flow, more than can be
controlled by a small number of functionals like drag and
lift. Therefore, an adaptive monitoring of the complete flow field, as
done by the multiscale analysis, is very desirable.

Here we develop our strategy for a test case. Since we focus on timestep
adaptation we will use uniformly refined meshes in space. 
Starting with a very coarse spatial mesh and CFL below unity, we
gradually establish sequences of timesteps which are well adapted to
the physical problem at hand. The scheme detects stationary regions,
where it switches to very high CFL numbers, but reduces the time
steps appropriately as soon as a perturbation enters the flow field.

Depending on the CFL number and the cost of the nonlinear solver, the adaptive
scheme chooses either explicit or implicit timesteps. For reasons of efficiency,
very small timesteps $\cfl \ll 1$ may be merged into a single step. This strategy
is detailed in Section \ref{sec.applications.strategy}.

Once we arrive at the fine spatial mesh, on which we really want to
compute and where most of the work is being done, we already work with
a very efficient time step. Moreover, we have a rational criterion
what the finest grid should be.

The paper is organized as follows: in Section \ref{sec.theory} we review the
theoretical background for our adaptive timestep control: DG and FV methods,
control of target functionals, error representation, space-time splitting, error
estimates. The new conservative approach for solving the dual problem is
presented in Section~\ref{sec.newapproach}. In Section \ref{sec.applications} we
define our adaptive strategy and apply it to compute perturbations of a
stationary shock. Some conclusions are drawn in Section \ref{sec.conclusions}.

%%%%%%%%%%%%%%%%%%%%%%%%%%%%%%%%%%%%%%%%%%%%%%%%%%%%%%%%%%%%%%%%%

\section{Derivation of space-time-split error estimates}
\label{sec.theory}
In this section, we recall some of the theoretical background of
adjoint error control, and we represent the extensions needed in our
time adaptive strategy. 
In Section~\ref{sec.dg_method} we introduce
the DG method used in the paper. In Section
\ref{sec.error_representations} we state the adjoint based error
representation for target functionals. In
Section~\ref{sec.space_time_splitting} we introduce a variant of the
projections in
space and time which lead to a splitting of the error representation.
One part decreases when the spatial grid is refined, and
the other part decays with the timestep.
The corresponding decay rates are a crucial ingredient of the time-adaptation strategy.  This
strategy and its application will be presented in
Section~\ref{sec.applications} below.

%%%%%%%%%%%%%%

\subsection{Discontinuous Galerkin methods for conservation laws}
\label{sec.dg_method}
Let $D$ be an open connected subset of $\mathbb R^d$, $d\geq 1$, let
$I := (0,T)$ be the time interval and let $\OT := D\times I$ be the
space-time domain, with boundary $\Gamma$ and outside unit normal
$\nu$. We consider the system
\begin{align}
\partial_t u +\nabla f(u) &=0 \quad \mbox{in } \OT, \\ 
\quad \fn(u) &= \gn \quad \mbox{on }\Gin,  \label{eq.cl}
\end{align}
where $u=(u_1,\ldots,u_m)^T :\OT \rightarrow \mathbb R^d$ is the vector of conservative variables
and $f(u) = (f_1(u), \dots, f_d(u))$ the flux matrix, with $f_i \in
C^1(\mathbb R^m, \mathbb R^m)$. The vector 
\begin{align*}
\fn(u) := (f(u),u)\cdot\nu
\end{align*}
is the space-time normal flux across the boundary, and for scalar
equations, the inflow boundary is given by
\begin{align*}
\Gin := \{(x,t) \in \Gamma \mid \frac{d}{du}((f(u),u)\cdot\nu  < 0 \}.
\end{align*}
Note that \eqref{eq.cl} includes initial data, since $D \times \{0\} \subset \Gin$, and
$\fn(u(x,0)) = u(x,0)$. For systems of conservation laws, the definition of in- and outflow
boundaries may be generalized via characteristic decompositions
\cite{HartmannHoustonJCP2002,SueliHouston2003}.

Let us define a partition of our time interval $I$ into subintervalls
$I_n = (t_{n-1}, t_n)$, where
\begin{align*}
0=t_0 < t_1 < \ldots < t_n < \ldots < t_N=T.
\end{align*}
Later on this partition will be defined automatically by the
adaptive algorithm.  Furthermore we define a regular polygonal
spatial grid $\mathcal {T}_D = \bigcup_j \{D_j\}$ such
that $\overline D = \bigcup_j \overline{D_j}$. We denote the 
corresponding space time prisms by
$$
\Ojn := D_j  \times I_n.
$$ For future reference, we denote the outward unit normal vector to
$\Ojn$ by $\nu_j^n$ or simply $\nu$. Thus we have constructed a
subdivision
\begin{align*}
\mathcal{T}_{\OT} = \bigcup_{j,n} \{\Ojn\}
\end{align*}
of the computational domain $\OT$. The spatial discretisation
$\mathcal T_D$ can change adaptively from timestep to timestep,
and for each fixed time interval $I_n$, the timestep is global
(i.e. it is the same for all spatial cells $D_j$).

\begin{remark}
We do not admit local timesteps, since we want to couple our time-adaptive strategy
to standard Runge-Kutta Finite Volume methods and Runge-Kutta
Discontinuous Galerkin methods.
\end{remark}
On this grid we define the following function spaces: First, let $S_h(\OT)$
be the mesh dependent broken space of discontinuous piecewise $H^1$
functions defined on $\mathcal{T}_{\OT}$,
\begin{align}
S_h(\OT) := \left\lbrace u \mid u_{\mid_{\Ojn}}\in H^1(\Ojn), \forall 
\Ojn \in \mathcal{T}_{\OT}\right\rbrace .
\end{align}
Furthermore we denote by $S_h^{s,r}(\OT)$ the (locally) finite dimensional
space consisting of discontinuous piecewise polynomial functions of
degree $s$ in space and $r$ in time defined on $\mathcal{T}_h$
\begin{align}
S_h^{s,r}(\OT) := \left\lbrace u_h \mid u_{h}(\cdot,t)\in 
P_{s}(D_j), \forall t\in {I_n},
u_{h}(x,\cdot)\in 
P_{r}(I_n), \forall x\in {D_j},
\forall D_j \times I_n\in \mathcal{T}_{\OT}\right\rbrace ,
\end{align}
where $P_{r}(I_n)$ denotes the space of polynomials of degree $r$ on
$I_n$ and $P_{s}(D_j)$ the space of polynomials of degree $s$ on
$D_j$.  Given a cell $\Ojn$ and a point $(x,t)\in\Gjn$, we define the
inner ($u^+$) and outer ($u^-$) values of a function $u \in S_h(\OT)$ with
via
\begin{align}\label{eq.fnpm}
u^\pm(x,t) := \lim _{{\delta \searrow 0 +} \atop {\delta \nearrow 0 -}}
u((x,t)-\delta \nu_j^n).
\end{align}
Defining the DG method for nonlinear conservation laws, whose
solutions in general contain shock waves, requires a careful
application of the theory of weak solutions, which states that
for a weak solution $u$ and a continuously differentiable test
functions $v$,
\begin{align*}
-(u,\partial_t v)_{\Ojn} -(f(u),\nabla v)_{\Ojn} + 
(f_\nu(u),v)_{\Gjn} = 0 \quad \forall j,n.
\end{align*}
Thus we have to define the normal flux $\fn(u)$ at the cell
boundaries, where the approximate solution $u_h$ is
discontinuous. This can be done with the help of numerical flux
functions, which we denote by $\fns$. So suppose that
$(x,t)\in\Gjn\setminus\Gamma$ is contained in an interior edge. If
$(x,t) \in \partial D_j \times I_n$, so that the normal points into
the spatial direction, then the canonical choice for $\fns$ is an
approximate Riemann solver
\begin{align}\label{eq.fn1}
\fns := f(u_h^+,u_h^-,n_j),
\end{align}
where $n_j$ is the outer normal to $D_j$ (i.e. $\nu_j^n=(n_j,0)$).
We require that the flux $\fns$ is consistent and conservative
in the sense of Lax.
If, on the other hand, $(x,t)\in D_j \times \partial I_n$, so that the
normal points into the time direction and $\fn(u)=u$, then we simply
require that $\fns$ be a convex combination of $u_h(x,t^\pm)$. More
specifically, suppose that $t=t_n$. Then we set
\begin{align}\label{eq.fn2}
\fns := u_h^\ast(x,t_n) := (1-\theta) u_h(x,t_n^+) + \theta u_h(x,t_n^-) 
\end{align}
for some value $\theta \in [0,1]$.
Different values of $\theta$ will yield different time discretisations, e.g.
explicit Euler for $\theta =0$, implicit Euler for $\theta =1$, if we work
with piecewise constant ansatz functions.

On the boundary of the domain, i.e. for $(x,t)\in\Gamma$, we set
\begin{align}\label{eq.fn3}
\fns :=  \left\{ \begin{array}{ll}
\gn & \mathrm{if} \; (x,t)\in\Gin \\
\fnpu & \mathrm{if} \; (x,t)\in\Gout
\end{array} \right.
\end{align}

In the following definition we simply state the
resulting DG(s,r) method, which is a discontinuous method both in
space and time. This definition is very similar to, see e.g.
\cite{BarthLarson2002,CockburnHouShu1990,HartmannHoustonSISC2002,
SueliHouston2003} and the references therein.
\begin{definition}
(i) The abstract semilinear form 
$\mathcal N : S_h(\OT) \times S_h(\OT) \to \mathbb R$ is given by
\begin{align}\label{eq.dg_form}
\mathcal N \left(u_h, v_h\right)
:= \sum_{j,n}
\left\lbrace 
 (\partial_t u_h + \nabla f(u_h),v_h)_{\Ojn} 
 + (\fns - \fn(u_h^+),v_h^+)_{\partial\Ojn}
\right\rbrace.
\end{align}
(ii) Now the DG(s,r) finite element method for the system of
hyperbolic conservation laws \eqref{eq.cl} is defined as follows: Find
$u_h \in S_h^{s,r}(\OT)$, such that
\begin{align}
\mathcal N \left(u_h, v_h\right) =0  \quad \forall v_h \in S_h^{s,r}(\OT).\label{eq.dg}
\end{align}
\end{definition}
As usual, the variational formulation \eqref{eq.dg_form},
\eqref{eq.dg} can be exploited as follows: Given $u_h \in S_h(\OT)$,
$\mathcal N (u_h,\cdot)$ is a linear functional on $S_h(\OT)$. Thus it can
be represented by an element of $S_h(\OT)$, which we call $R(u_h)$, the
{\em residual}. On the interior of a cell $\Ojn$ we introduce the 
{\em cell residual} 
\begin{equation}\label{eq.r1}
R_h := \partial_t u_h + \nabla f(u_h)
\end{equation}
and on the boundaries $\Gjn$ the {\em edge residual}
\begin{equation}\label{eq.r2}
r_h := \fns - \fn(u_h^+).
\end{equation}
Then \eqref{eq.dg_form} can be rewritten as 
\begin{equation}\label{eq.r3}
(R(u_h),v_h) 
= \sum_{j,n}
\left\lbrace 
(R_h,v_h)_{\Ojn} + (r_h,v_h^+)_{\partial\Ojn}
\right\rbrace.
\end{equation}
The $DG(s,r)$  solution $u_h\in S_h^{s,r}(\OT)$ of \eqref{eq.dg}
is now given by
\begin{equation}\label{eq.galerkin_orthogonality}
(R(u_h),v_h) = 0 \quad \forall v_h \in S_h^{s,r}(\OT),
\end{equation}
which is the classical Galerkin orthogonality: the residual $R(u_h)$ is
orthogonal to the test space $S_h^{s,r}(\OT)$.

In the following, we mostly work with the $DG(0,0)$ and $DG(1,1)$ methods, both
in their explicit ($\theta=0$) and implicit ($\theta=1$) form.  The DG(0,0)
method is equivalent to a first order accurate finite volume scheme, using
explicit order implicit Euler scheme for the time integration. 
\iffalse
For higher oder
accurate finite volume schemes we refer to Barth and Larson 
\cite{BarthLarson2002}.  
\fi
In \cite{BarthLarson2002}, Barth and Larson derive a weak formulation of the
form \eqref{eq.galerkin_orthogonality} for higher oder accurate finite volume schemes.

Therefore, the techniques presented in this paper can be
applied to finite volume and Discontinuous Galerkin methods.

\subsection{Adjoint error representation for target functional}
\label{sec.error_representations}
In this section we define the class of target functionals treated in
this paper, state the corresponding adjoint problem and recall the
classical error representation which we will later use for
adaptive time step control.

Our objective is to estimate the error in a user specified functional
$J(u)$, which can be expressed as a sum of weighted integrals over the
domain $\OT$ and the outflow boundary $\Gout$.
Typical examples of such functionals are the lift or the drag of a
body immersed into a fluid. 

To simplify matters we consinder functionals of the following form:
\begin{align*}
J(u) = (u,\psi)_{\OT} - (\fn(u),\psi_\Gamma)_{\Gout}
\end{align*}
Our purpose is to control the error
\begin{align*}
J(u)-J(u_h) .
\end{align*}
In order to derive the classical error representation one linearizes
the evolution equation satisfied by the error $u-u_h$ and works with
the adjoint equation of the linearized error equation. Thus we
introduce an approximate Jacobian $\ab$ of $f$ by
\begin{align}
\ab &:= \int\limits_0^1 \frac{d}{d\tau}f(u_h+\tau(u-u_h))d\tau. \label{eq.ab}
\end{align} 
Note that
\begin{align*}
f(u)-f(u_h) &= \ab (u-u_h). 
\end{align*}

In practice we linearize around the approximate solution.
A direct calculation yields the following theorem:
\begin{theorem}\label{theorem.error_representation}
Suppose $\varphi \in H^1(\Omega)$ solves the adjoint problem
\begin{align}\label{eq.adjoint_1}
\varphi &= \psi_\Gamma \quad \mathrm{on} \; \Gout \\
\label{eq.adjoint_2}
\partial_t \varphi + \ab \nabla \varphi &= \psi \quad \; \; \mathrm{in} \; \Omega.
\end{align}
Then for all $\varphi_h \in S_h^{s,r}(\OT)$, the error in the target functional satisfies 
\begin{align}\label{eq.error_representation}
J(u)-J(u_h) = ( R(u_h),\varphi-\varphi_h).
\end{align}
\end{theorem}
Equivalently one can also define the adjoint solution via a
variational formulation (see e.g.
\cite{BarthLarson2002,HartmannHoustonSISC2002,SueliHouston2003}).
In \cite{Tadmor1991} Tadmor proves the well-posedness of the
adjoint problem \eqref{eq.adjoint_1} -- \eqref{eq.adjoint_2}
for scalar, convex, one-dimensional conservation laws. 
The key observation is that, if the forward solution $u$ has jump
discontinuities, then due to the entropy condition the jump of the
transport coefficient $\ab$ has a distinct sign. This makes it
possible to follow the characteristics of the adjoint problem
backwards in time.

Identity \eqref{eq.error_representation} is the error representation
which we discussed in the introduction and onto which we are going to
base our adaptive strategy. By definition \eqref{eq.r1} -
\eqref{eq.r3} of the residual $ R(u_h)$, the error
representation may be decomposed as a sum over the cells and edges of
inner products of the local residuals with the solution of our dual
problem. Due to Galerkin orthogonality
\eqref{eq.galerkin_orthogonality}, we can subtract an arbitrary
test function $\varphi_h$, which is very convenient when we derive local
error estimates later on.

%%%%%%%%%%%%%%%%%%%%%%%%%%%%%%%%%%%%%%%%%%%%%%%%%%%%%%%%%%%%%%%%%%%%%%%%%%%%%%%%
\subsection{Space-time splitting and the error estimate}
\label{sec.space_time_splitting}

The error representation \eqref{eq.error_representation} is not yet
suitable for time adaptivity, since it combines space and time
components of the residual and of the difference $\varphi-\varphi_h$ of the dual
solution and the test function. The main result of this section is an error estimate whose components
depend either on the spatial grid size $h$ or the time step $k$, but
never on both. The key ingredient is a space-time splitting of
\eqref{eq.error_representation} based on $L^2$ projections.
Similar space-time projections were introduced previously
in \cite{Hartmann1998, Sueli}.
Here we adapt them to the finite element spaces used in
the error representation \eqref{eq.error_representation}.

Let $P_{s,r}(\Ojn) = P_s(D_j)\times P_r(I_n)$ be the space of polynomials of
degree $s$ on ${D_j}$ and $r$ on $I_n$.
Furthermore let
$\hat P_{I_n}^r(\Ojn)=\{w\in L^2(\Ojn)|w(x,\cdot)\in P_r(I_n),\forall x\in D_j\}$, and
$\hat P_{D_j}^s=\{w\in L^2(\Ojn)|w(\cdot,t)\in P_s(D_j),\forall t\in I_n\}$.
For $r \geq 0$ define the $L^2$ projection 
$\Pi_{I_n}^r: L^2(\Ojn) \to \hat P_{I_n}^r(\Ojn)$ via
\begin{align}
(u-\Pi_{I_n}^r u, \varphi)_{I_n} = 0 \quad \forall \varphi\in \hat P_{I_n}^r(\Ojn), \forall x\in D_j,
\end{align}
and for $s \geq 0$ define the $L^2$ projection
$\Pi_{D_j}^s: L^2(\Ojn) \to \hat P_{D_j}^s(\Ojn)$ via
\begin{align}
(u-\Pi_{D_j}^s u, \varphi)_{D_j} = 0 \quad \forall \varphi\in \hat P_{D_j}^s(\Ojn), \forall t\in I_n.
\end{align}
Similarly let the $L^2$ projection
$\Pi_{\Ojn}^{s,r}: = L^2(\Ojn) \to  P_{s,r}(\Ojn)$
be defined via
\begin{align}
(u-\Pi_{\Ojn}^{s,r}u,\varphi)_{\Ojn}=0\quad \forall \varphi\in 
P_{s,r}(\Ojn).
\end{align}
Note that $\Pi_{\Ojn}^{s,r} = \Pi_{{D_j}}^s \Pi_{I_n}^r$.

First we choose $\varphi_h$ in the error representation
(\ref{eq.error_representation}) to be $\varphi_h =\Pi_{h,k}^{s,r} \varphi$,
i.e. $\varphi_h\mid_{\Ojn} =\Pi_{{D_j}}^s \Pi_{I_n}^r \varphi= \Pi_{I_n}^r \Pi_{{D_j}}^s
\varphi$, with $\Pi_{I_n}^r$ and $\Pi_{{D_j}}^s$ as defined above.
Using the identity
\begin{align*}
\varphi-\Pi_{h,k}^{s,r}\varphi =\varphi-\Pi_{I_n}^r \varphi + \Pi_{I_n}^r
\varphi -\Pi_{h,k}^{s,r}\varphi = (id-\Pi_{I_n}^r 
)\varphi+(id-\Pi_{D_j}^s )\Pi_{I_n}^r \varphi
\end{align*}
we obtain the following splitting of the error representation:
\begin{align}\label{eq.error_representation2}
J(u)-J(u_h) &= ( R(u_h),(id-\Pi_{I_n}^r 
)\varphi+(id-\Pi_{D_j}^s )\Pi_{I_n}^r \varphi)\\ 
&= \sum_{j,n}
\{ \underbrace{(R_h,(id-\Pi_{I_n}^r)\varphi)_{\Ojn} +
(r_h,(id-\Pi_{I_n}^r)\varphi^+)_{\partial\Ojn}}_{\eta_{k}^{jn}}\\
&\phantom{ = \sum_{j,n}}+\underbrace{(R_h,(id-\Pi_{D_j}^s )\Pi_{I_n}^r \varphi)_{\Ojn} +
(r_h,(id-\Pi_{D_j}^s )\Pi_{I_n}^r \varphi^+)_{\partial\Ojn}
}_{\eta_{h}^{jn}}\}\\ \label{eq.error_representation3}
&=: \eta_k +\eta_h,
\end{align}
where $\eta_k$ is the time-component and $\eta_h$ the space-component of the
error representation $\eta$.

In this paper we consider grids, which are locally tensor products of a
spatial grid $\mathcal {T}_{D_j}$ and a timestep $I_n$. For an implicit Runge-Kutta Finite Volume Method
$\eta_{k}^{jn}$ and $\eta_{h}^{jn}$ then take the form

\begin{align*}
\eta_{k}^{jn}=&  
 (R_h, (id-\Pi_{I_n}^r)\varphi)_{D_j \times I_n} 
 +(\fns - \fn(u_h^+),(id-\Pi_{I_n}^r)\varphi^+)_{\partial D_j \times I_n} \\
&+ (\left[ u_h \right]_{n-1} ,(id-\Pi_{I_n}^r)\varphi_{n-1}^+)_{D_j}\\
\eta_{h}^{jn} =& (R_h, (id-\Pi_{D_j}^s )\Pi_{I_n}^r \varphi)_{D_j \times I_n} 
+(\fns - \fn(u_h^+),((id-\Pi_{D_j}^s )\Pi_{I_n}^r \varphi)^+)_{\partial D_j \times 
 I_n} \\
&+ (\left[ u_h \right]_{n-1} ,((id-\Pi_{D_j}^s )\Pi_{I_n}^r 
\varphi)_{n-1}^+)_{D_j}
\end{align*}

where the flux difference on the spatial boundaries $\partial D_j \times I_n$
and the jump of $u_h$ on the time boundary $D_j\times\{t_{n-1}\}$ are
realizations of the residual term $r_h$ in \eqref{eq.r2}.

For future reference, we also introduce the quantities
\begin{align}
\eta_k &:=  \sum_{j,n}\eta_{k}^{jn}&
\eta_h &:=  \sum_{j,n}\eta_{h}^{jn}\\
\label{eq.bar_eta_kn}
\bar\eta_k^n&:= \frac{1}{k_n}\sum_{{D_j}\in \mathcal T}|\bar\eta_{k}^{jn}|&
\bar\eta_h^n&:= \frac{1}{k_n}\sum_{{D_j}\in \mathcal T}|\bar\eta_{h}^{jn}|\\
\bar\eta_k &:=  \sum_{n}k_n|\bar\eta_{k}^{n}|&
\bar\eta_h &:=  \sum_{n}k_n|\bar\eta_{h}^{n}| \\
\bar\eta &:=\bar\eta_k +\bar\eta_h.&
\end{align}

In Section \ref{sec.SEP} we will show numerically, that the error terms
$\bar\eta_{k}$ and $\bar\eta_{h}$ depend on $k$ and $h$.

\section{A new approach to solving the adjoint problem}
\label{sec.newapproach}

The error representation \eqref{eq.error_representation} assumes that the exact
solution $\varphi$ of the dual problem \eqref{eq.adjoint_2} is available. This
is, of course, not the case. All we can do is to compute an approximation 
$\vps$ of $\varphi$. An important question is in which space we should choose
the approximation $\vps$ (let us call this space $S_\sharp$ ). If we choose $S_\sharp
\subseteq S_h^{s,r}$, then - due to Galerkin orthogonality of the residual - the
error representation \eqref{eq.error_representation} would return zero.
Therefore, $S_\sharp$ should not be contained in $S_h^{s,r}$.

There are essentially three approaches in the literature to compute an
approximate solution to the dual problem. The first approach is to keep the
polynomial degrees $r$ and $s$ fixed, but compute the solution to the dual
problem on a  finer grid $\mathcal {T}_{D_j} \subset \mathcal {T}_{D _{j+1}} $.
The second approach is to compute the dual solution using higher order finite
elements and using projections to get $\varphi_h$: 
\begin{align*}
\mathrm{Compute\!:} \;\; \vps \in S_h^{s+1,r+1}(\OT)
\qquad\rightsquigarrow\qquad
\varphi_h:= \Pi \vps, 
\end{align*}
where $\Pi$ is the projection from the higher order finite element space onto the test space $S_h^{s,r}(\OT)$.
The third way is to compute a solution in the test space of the forward problem, which means to use
the same order finite elements,  and then do a higher order reconstruction $R$.
\begin{align*}
\mathrm{Compute\!:}\;\;\varphi_h \in S_h^{s,r}(\OT)
\qquad\rightsquigarrow\qquad
\vps := R\varphi_h
\end{align*}

In the following, we describe a fourth approach, which avoids to approximate
$\varphi$ alltogether. Instead, we approximate the spatial gradient
$\nabla_x\varphi$. The remarkable fact is that this gradient satisfies a
conservation law instead of a nonlinear transport equation, and its numerical
approximation is therefore very robust in the presence of shocks. 
In the present paper we limit our presentation to first order schemes in one
space dimension. Our approach can be applied to the dual problem, if the forward
problem is approximated by a first order DG method, or a Finite Volume method.
The backward problem can then be computed by the same method as the forward
problems. The generalization of our ansatz to higher order schemes is relatively
straightforward in one space dimension.

Let us look at the details: Due to Galerkin orthogonality, the dual solution
$\varphi$ does not enter the error representation as such. Instead, the relevant
term is the difference of the dual solution and its projection to the finite
element space, $\varphi-\varphi_h$.  Using one of the three methods described
above, one needs additional degrees of freedom to compute an approximation
$\varphi$ to the dual problem, and some computed information will never be used,
since only the difference  $\varphi -\varphi_h$ enters the error representation.
Therefore we suggest to compute the spatial gradient of the dual solution. 

To illustrate our approach (still in one spatial dimension), we assume that
$\varphi_h$ is the piecewise constant function satisfying
\begin{eqnarray*}
\varphi_h(x,t) &\equiv& \varphi(x_0, t_0) \quad\mathrm{for}\quad (x,t)\in
D_j\times I_n.
\end{eqnarray*}
for some given point $(x_0, t_0)\in D_j\times I_n$ (e.g. the midpoint).
Expanding $\varphi$ around $(x_0,t_0)$,
\begin{eqnarray*}
\varphi(x,t) &=&\varphi(x_0,t_0)+(x-x_0) \partial_x \varphi(x_0, t_0) +(t-t_0)
\partial_t
\varphi(x_0, t_0)+O(h^{2}+k^2),
\end{eqnarray*}
and using the adjoint equation \eqref{eq.adjoint_2}, we obtain that
\begin{align*}
\varphi-\varphi_h(x,t) &=(x-x_0) \partial_x \varphi(x_0, t_0) +(t-t_0) 
	\partial_t \varphi(x_0, t_0)+O(h^{2}+k^2)\\
      &=(x-x_0) \partial_x \varphi(x_0, t_0)+(t-t_0) 
      (\psi-\ab \partial_x \varphi(x_0,t_0 )) + O(h^{2}+k^2)\\
      &=[(x-x_0)+(t-t_0)(\psi-\ab)] \partial_x \varphi(x_0, t_0)
       + O(h^{2}+k^2).
\end{align*}
Since $\psi$ and $\ab$ are assumed to be known, the only unknown function is $\partial_x
\varphi(x_0, t_0)$.

In order to derive the differential equation which is satisfied by 
$\partial_x\varphi$, we differentiate the adjoint equation \eqref{eq.adjoint_2},
\begin{align*}
\varphi &= \psi_\Gamma \quad \mathrm{on} \; \Gout \\
\partial_t \varphi + \ab \partial_x \varphi &= \psi \quad \; \; \mathrm{in} \; \Omega
\end{align*}
with respect to $x$ and obtain
\begin{align}
w &= \partial_x \psi_{\mathit{\Gamma}}   \quad  \mathrm{on} \; \Gout \label{eq.dp2.1}\\
\partial_t  w + \partial_x (\ab w) &= \partial_x \psi  \quad \; \; \mathrm{in} \;
\Omega,
\label{eq.dp2.2}
\end{align}
where $w:= \partial_x \varphi$.

Therefore it is not necessary to compute the approximations
$\vps$ and $\varphi_h$ of $\varphi$, but it is sufficient to compute an
approximation $w_\sharp \in S_h^{s,r}(\OT)$ of $\partial_x \varphi$.

\begin{remark}
It is striking to note that the gradient $w= \partial_x \varphi$
actually satisfies a conservation law, \eqref{eq.dp2.1}-\eqref{eq.dp2.2},
instead of a linear transport equation,
\eqref{eq.adjoint_1}-\eqref{eq.adjoint_2}. Therefore, $w_\sharp$ can be
computed  with the same algorithm as the forward problem, and in the same finite
element space. This leads to an efficient and robust solver: for discontinuous
$\ab$, finite difference schemes for \eqref{eq.adjoint_1}-\eqref{eq.adjoint_2}
may suffer from serious stability problems. Due
to their upwind nature, finite volume schemes for the conservation law
\eqref{eq.dp2.1}-\eqref{eq.dp2.2} handle  discontinuous coefficients easily.

In work in progress, we are analysing the efficiency of the new approach in more
detail, generalize it to higher order and several space dimensions, and study
related issues like boundary conditions for compressible fluid flows.
\end{remark}

We will use this new approach in the numerical examples in Section
\ref{sec.applications}.

%%%%%%%%%%%%%%%%%%%%%%%%%%%%%%%%%%%%%%%%%%%%%%%%%%%%%%%%%%%%%%%%%%%%%%%%%%%%

\section{Time adaptive strategy and application to perturbed shocks}
\label{sec.applications}

In this section we describe the strategy for adaptive time step
control, define a suitable numerical experiment and present first
numerical results which demonstrate the potential of this approach.

\subsection{The adaptive strategy}
\label{sec.applications.strategy}

In many applications, there are canonical target functionals which are of great interest to the user, like the lift and drag in aerodynamics. In some cases, an error margin may be prescribed for a given application. In other cases, it is less clear which accuracy should be and can be provided by a numerical computation, and with reasonable resources. In the following, we suggest prototype strategies to deal with both situations, where the tolerance may be, or may not be, prescribed. Many equally valid variants of these could be proposed, as well. As pointed out before, we focus on the time adaptation. For clarity of exposition, we therefore use uniformly refined spatial grids.

In the present paper, we only treat Burgers' equation. In a paper in preparation, we extend this to the
Euler equations of gas dynamics. We begin by computing the forward and the dual solution as well as the
error estimator on a relatively coarse spatial grid ($L=0$). Usually this spatial grid is much coarser than
the grid we actually want to compute on. Since we want to compute a solution with accuracy comparable to an explicit solution, we prescribe a uniform CFL number below unity in this first computation (e.g. CFL=0.8).

After evaluating the error representation, we have to take two decisions:

\begin{enumerate}
\item
the refinement level $L$ of the next spatial grid. In some cases we will gradually increase the level by 
one. This careful approach may be important if it is not clear whether the dynamics of the solution is already captured on the present grid. In other cases (including the example treated below), the time dynamics is already resolved very well on level $L=0$, and we can immediately proceed to the finest grid level.
\item
the tolerance $\tol_k(L)$ for the temporal component of the error,
$\bar\eta_k$. The choice of $\tol_k(L)$ will be based on assumptions of the
asymptotic decay of the error. If, as in Figure \ref{fig.burgers_1_etabar}, 
the error decays to first order, then we may choose $\tol_k(L+1)=0.5\,\tol_k(L)$.
\end{enumerate}

%\begin{align*}
%|\bar{\eta_k}| < Tol \, |J(u)|.
%\end{align*}

Now we adapt the timestep locally in order to equidistribute the error
densities $\bar\eta_k^n$. Recall from \eqref{eq.bar_eta_kn} that 
$$
	\bar\eta_k^n = \frac{1}{k_n}\sum_{{D_j}\in \mathcal T}|\bar\eta_{k}^{jn}|,
$$
and
$$
	\sum\limits_n k_n \bar\eta_k^n = \bar\eta_k.
$$
If the $\bar\eta_k^n$ were already equidistributed with respect to $n$, then
they would satisfy
$$
	%\bar\eta_k^n = \frac1T \, \bar\eta_k \quad\mathrm{for \; all} \; n.
	\bar\eta_k^n = \bar\eta_k/T \quad\mathrm{for \; all} \; n.
$$
Now, instead of aiming at local error densities of $\bar\eta_k/T$, we
target at an equidistribution of
$$
	\bar\eta_k^n \approx \tol_k(L+1)/T \quad\mathrm{for \; all} \;
	n,
$$
where $\tol_k(L+1)$ is a given tolerance on grid $(L+1)$.
Assuming once more that the time component of the error varies
linearly with the time step, we compute the new timestep $k_m$ (on level
$(L+1)$) as 
\begin{align}
	\label{eq.k_m}
	k_m :=k_n \; \frac{
	\tol_k(L+1)/T}{\bar\eta_k^n} \; .
\end{align}
Using this new timestep distribution we perform a new computation
on the finer spatial grid. Note that due to the linear decay of the error with the timestep, often the new
distribution has a similar number of timesteps as the previous one.

If a total tolerance for the error, $|J(u)-J(u_h)|<\tol_{tot}$ is prescribed, 
then the above loop is stopped once
$$
	\bar\eta_k+\bar\eta_h< \tol_{tot}.
$$
Our experience so far is the
following: already on very coarse grids, the method detects the areas of
stationary and instationary flow quite well, and chooses the time steps
accordingly.

We would like to call the approach which combines \eqref{eq.k_m} with an 
implicit solver the {\em adaptive, fully implicit strategy}. A possible drawback
of this strategy is that it may lead to extremely small timesteps (CFL $\ll 1$)
when strong instationary waves pass the computational domain. Therefore, in the
second and third example, we restrict the time step size from below. When the
equidistribution of the error suggests $\cfl<5$, we switch to an explicit solver
with $\cfl =0.8$. This saves a considerable number of timesteps. We call this
approach the {\em adaptive, implicit/explicit strategy}.

%%%%%%%%%%%

\subsection{Test problem and asymptotic decay rates}
\label{sec.applications.test_problem}

%{Beispiel 4, timeit(48,0.7759,0), s11}

Now we set up an instationary test case, which is almost stationary, such that
an implicit (or implicit/explicit) scheme might be superior to a fully
explicit one. Our choice is a perturbed stationary shock for Burgers' equation
\begin{align*}
u_t +(\oh u^2)_x =0
\quad \mathrm{for} \; x\in[0,1] \; \mathrm{and} \; t\in[0,48].
\end{align*}
The initial data and corresponding unperturbed solution are given by
\begin{align*}
u(x,t) = 
\left\{ \begin{array}{rrl}
 1 & \mathrm{for} & x<0.5 \\
-1 & \mathrm{for} & x>0.5
\end{array}\right. .
\end{align*}
Then we place a disturbance at the left boundary of the
domain, which makes the stationary shock move. The new 
shock position as a function of time is given by
\begin{align*}
s(t) = 
\left\{ \begin{array}{rrl}
 0.5 & \mathrm{for} & t<12 \\
 0.5+\theta_1(t) sin(\frac{2\pi}{3}(t-12)) & \mathrm{for} & 12<t<18 \\
 0.5 & \mathrm{for} & 18<t<30 \\
 0.5+\theta_2(t) sin(\frac{2\pi}{3}(t-30)) & \mathrm{for} & 30<t<36 \\
 0.5 & \mathrm{for} & 36<t 
\end{array}\right..
\end{align*}
where 
\begin{align*}
\theta_1(t) = 7.5 \cdot 10^{-3} (t-12)^4(t-18)^4 / {6561} \\
\theta_2(t) = 0.5 \cdot 10^{-3} (t-30)^4(t-36)^4 / {6561}
\end{align*}

\begin{figure}[ht]
\begin{center} 
\epsfig{file=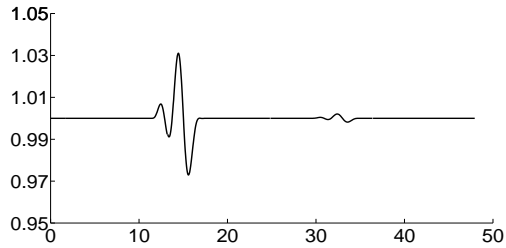, width=0.49\linewidth,height=0.15\textheight}
\end{center}
\caption{Burgers' equation: Left boundary data for perturbed stationary shock.} 
\label{fig.burgers_bsp3}
\end{figure}

Using characteristic theory, we can derive the perturbed left boundary
condition, which is displayed in Figure~\ref{fig.burgers_bsp3}. Note that
the magnitude of the first perturbation is about 1.5 percent of the
shock strength and that of the second perturbation about 0.1 percent.
The functional $J(u)$ is a weighted mean value in space
and time of the solution,
\begin{align*}
J(u) := \int_0^T\int_{0.25}^{0.65}u(x,t)\exp\left(-\frac{1}{1-y(x)^2}\right)dxdt,
\end{align*}
where $y(x):= (x-0.45)/0.2$. Note that the integration area completely
covers the domain containing the shock.

\subsubsection{Asymptotic decay rates}%{Separate refinement in space and time}
\label{sec.SEP}
Since the adaptive strategy outlined in Section \ref{sec.applications.strategy}
above depends on assumptions on the assymptotic behavior of the error, we first 
try to estimate these decay rates. There is no analytical result which shows how the error terms
$\bar\eta_k$ and $\bar\eta_h$ depend on $k$ and $h$. Therefore, we estimate this 
dependence numerically. We compute the perturbed shock described in
Section~\ref{sec.applications.test_problem} with a first order
finite volume method with Engquist-Osher flux, which is equal to a DG(0,0) method.
We compare the two approaches:
\begin{itemize}
\item refinement only time 
\item and refinement only space.
\end{itemize}

\begin{figure}[hbtp]
\begin{center} {
 \subfigure[uniform refinement in time]
 {\epsfig{file= 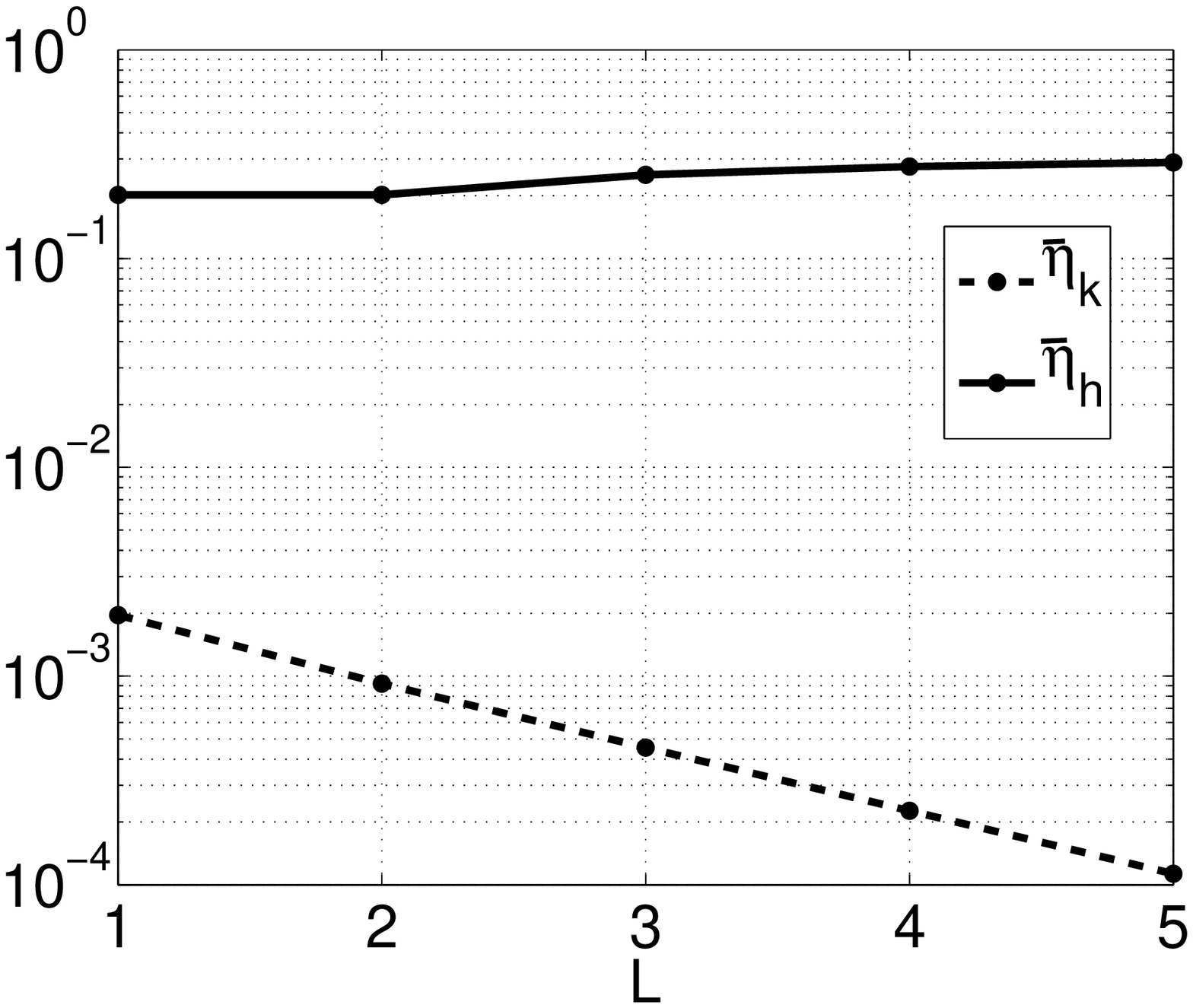, width=0.49\linewidth}}
 \subfigure[uniform refinement in space]{\epsfig{file= 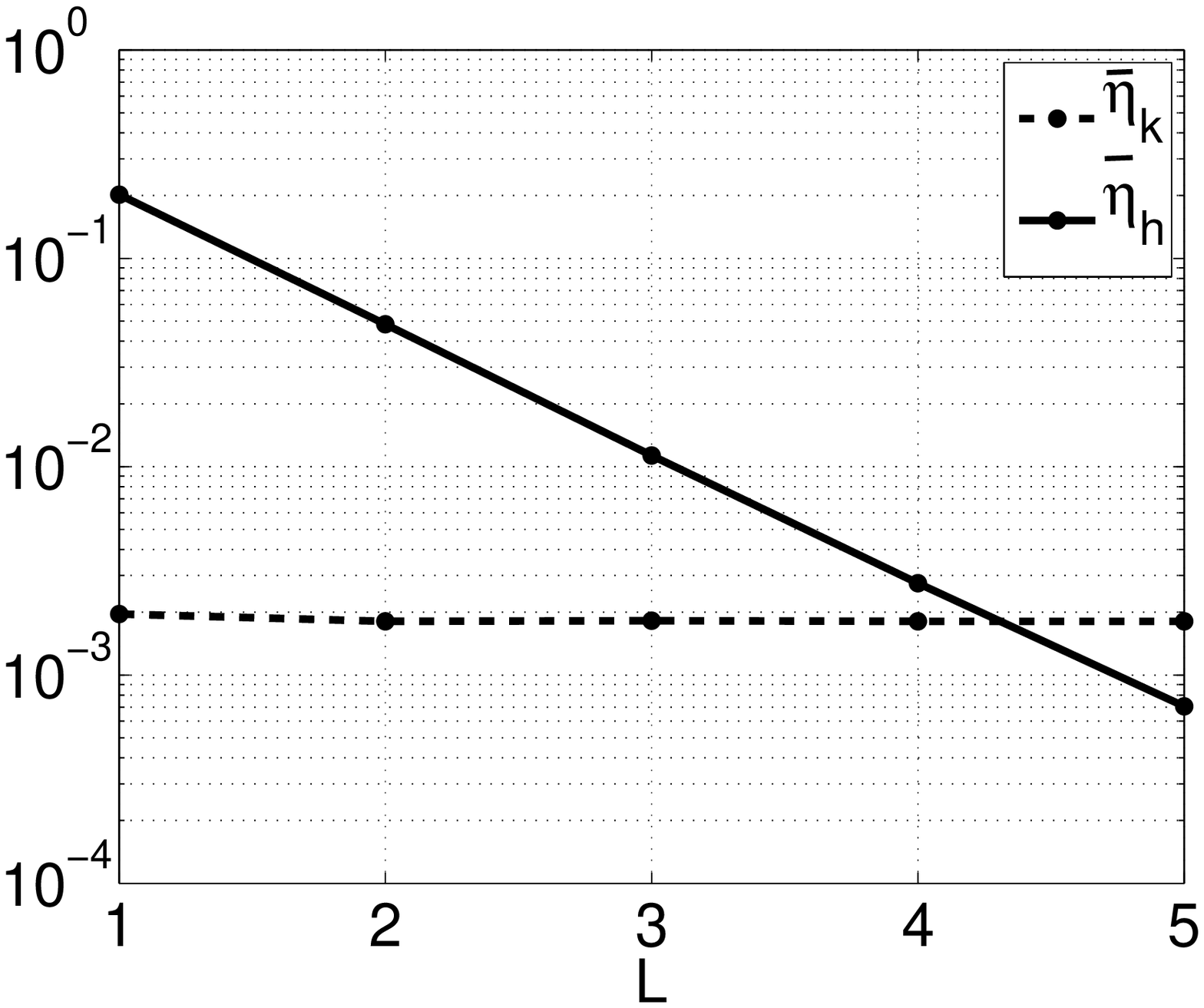, width=0.49\linewidth}}
}

\end{center}
\caption{Error representation for Burgers equation, 
first order method, $\bar{\eta_k}$ and $\bar{\eta_h}$ versus level 
 of refinement. (a) uniform refinement in time. (b) uniform refinement in space.} 
 \label{fig.burgers_1_etabar}
\end{figure}

\iffalse
\begin{figure}[hbtp]
\begin{center} {
 \subfigure[uniform refinement in time, level of refinemet]{\epsfig{file= pic3/ref_time2.eps, width=0.49\linewidth}}
 \subfigure[uniform refinement in space, level of refinemet]{\epsfig{file= pic3/ref_space2.eps, width=0.49\linewidth}}
}
\end{center}
\caption{Error representation for Burgers equation, 
first order method, ${\eta_k}$, ${\eta_h}$.} \label{fig.burgers_1_eta}
\end{figure}
\fi

\begin{table}  
\begin{tabular}{c|c|c|c|c|c|c|c|c|c}
$L$ &$dx$ & $dt$ & $\bar{\eta_k}$	&$\bar{\eta_h}$&$\eta_k$	&$\eta_h$	&$
J(u_h)$ &$\eta_h+\eta_k$& $\theta$ \\ \hline
1&  0.050000&  0.038795&    1.96e-03&  2.02e-01&  1.29e-04&  2.00e-01&  1.72e+00&  2.01e-01&  5.57e+00   \\  \hline 
2&  0.025000&  0.019398&    9.81e-04&  4.83e-02&  1.83e-05&  4.75e-02&  1.74e+00&  4.75e-02&  5.49e+00   \\  \hline 
3&  0.012500&  0.009699&    4.81e-04&  1.21e-02&  3.57e-06&  1.17e-02&  1.75e+00&  1.17e-02&  5.71e+00   \\  \hline 
4&  0.062550&  0.004849&    2.37e-04&  3.10e-03&  1.30e-06&  2.89e-03&  1.75e+00&  2.89e-03&  7.06e+00
\end{tabular}\\ 
\caption{Efficiency $\theta = \frac{\eta_h+\eta_k}{J(u)-J(u_h)}$of the error representation} 
\end{table}

Each of the plots in Figure~\ref{fig.burgers_1_etabar} show the error
estimators $\bar\eta_k$ (error in time) and  $\bar\eta_h$ (error in space).  In the
Figure~\ref{fig.burgers_1_etabar}(a) we refined only in time. Here the spatial
error remains constant, while the time error still decreases with first order.
The second Figure \ref{fig.burgers_1_etabar}(b) shows the refinement only in space. 
The time error $\bar\eta_k$ is almost constant, while the spatial error is decreasing with second order. 

Numerically the terms $\bar\eta_t$  and $\bar\eta_h$ behave as expected. 
They depend either on $k$ or on $h$, but never on both. The behaviour of
$\eta_h$ and $\eta_k$ is very similar, and not displayed here.

\begin{remark}
The numerically validated results can be used for adaptive 
grid refinement. The error estimator $\bar\eta_h$ can be used as an indicator for spatial
adaption and the estimator $\bar\eta_k$ for time step control.
\end{remark}

\subsection{Computational results}
\label{sec.applications.results}

{\bf Example 1:}   The first computation ($L=0$) is done on a grid with 20
spatial cells and a uniform CFL number of 0.8 using explicit timesteps. It needs
$N=1238$ timesteps, reaching a total error of $\bar\eta =0.204$ and a relative
error of $|\bar\eta/J(u_h)| = 11.9 \%$,  but a temporal error of
$|\bar\eta_k/J(u_h)| = 0.13 \%$. Our adaptive strategy now aims at a time step
distribution on the next grid with tolerance $\tol_k(L+1) = \bar\eta_k(L)$.
Based on the assumption that the time component of the error varies linearly
with the time step (which is motivated by  Fig.~\ref{fig.burgers_1_etabar}), the
scheme chooses new timesteps on the next grid according to the equidistribution
rule \eqref{eq.k_m}. 

The second row of Table~\ref{table.perturbed_shock.1}, for level $L=1$, gives
also $N=1238$ time steps, now using adaptive implicit timesteps. Now the
relative temporal error is $|\bar\eta_k/J(u_h)| = 0.073 \%$, and it is dominated
by the spatial error $|\bar\eta_h/J(u_h)| =  2.2 \%$.

\begin{table}  
\begin{tabular}{c|c|c |c|c} 
L	& N	&$\bar\eta_k/J(u_h)$	&$\bar\eta_h/J(u_h)$	&$\bar\eta/J(u_h)$	\\ \hline  	 
0	&1238	&1.34e-03		&1.17e-01		&1.19e-01	\\  
1	&1238	&7.37e-04		&2.20e-02		&2.27e-02	\\  
\end{tabular}\\  
\caption{Example 1: Perturbed shock for Burgers' equation.  
From left to right: level $L$, number of time steps $N$, time component
of error estimator $\bar\eta_k/J(u_h)$, spatial component
of error estimator $\bar\eta_h/J(u_h)$, total error estimator
$\bar\eta/J(u_h)$.}
\label{table.perturbed_shock.1}
 
%\vspace*{-4ex}

\end{table}

%%%%%%%%%%

\begin{figure}[hbtp]
\begin{center} {
 \subfigure[$\cfl(t_n)$, L=0, uniform explicit timestep]
 	{\epsfig{file= 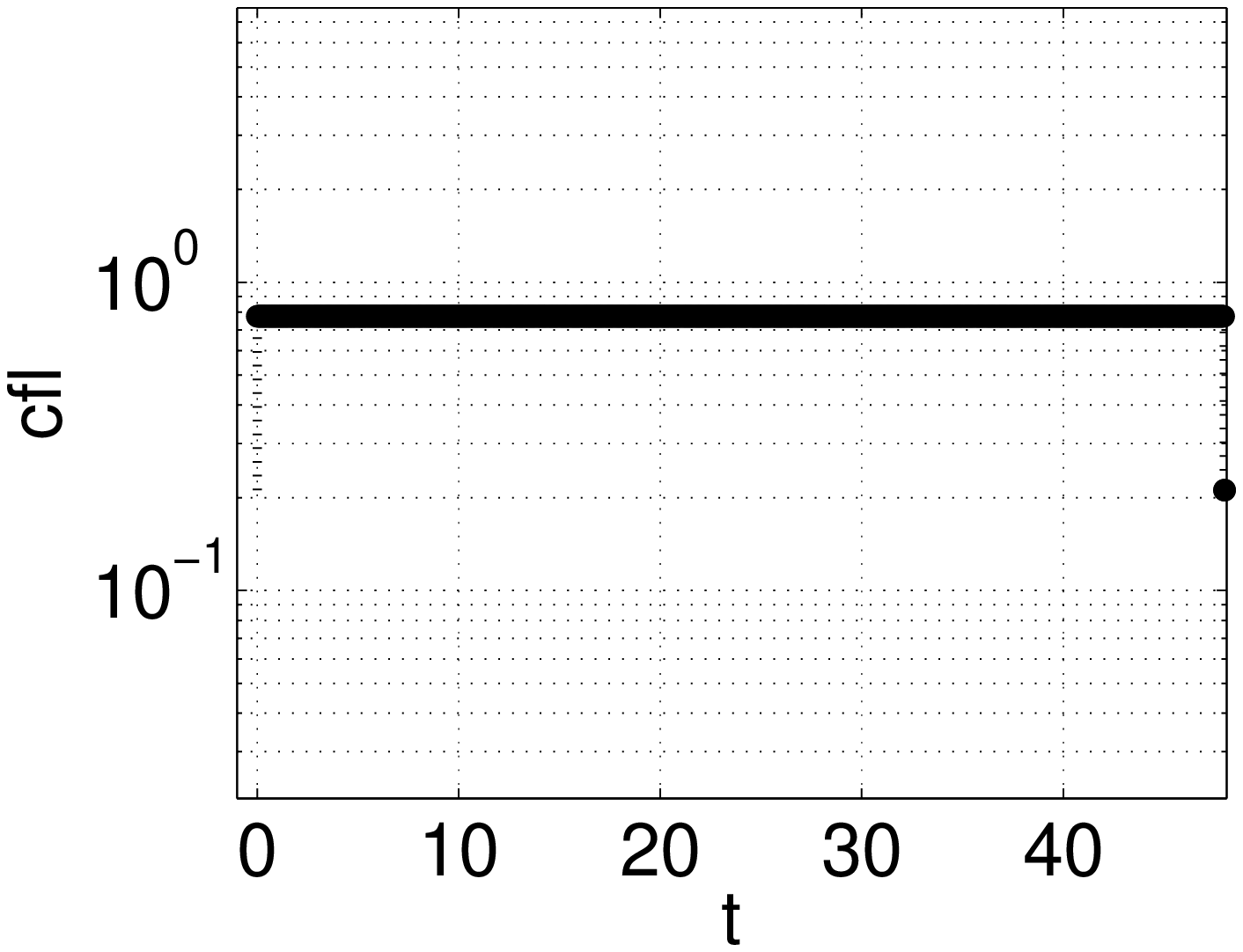, width=0.49\linewidth,height=0.16\textheight}}
 \subfigure[$\bar\eta_k^n(t_n)$, L=0, uniform explicit timestep]
 	{\epsfig{file= 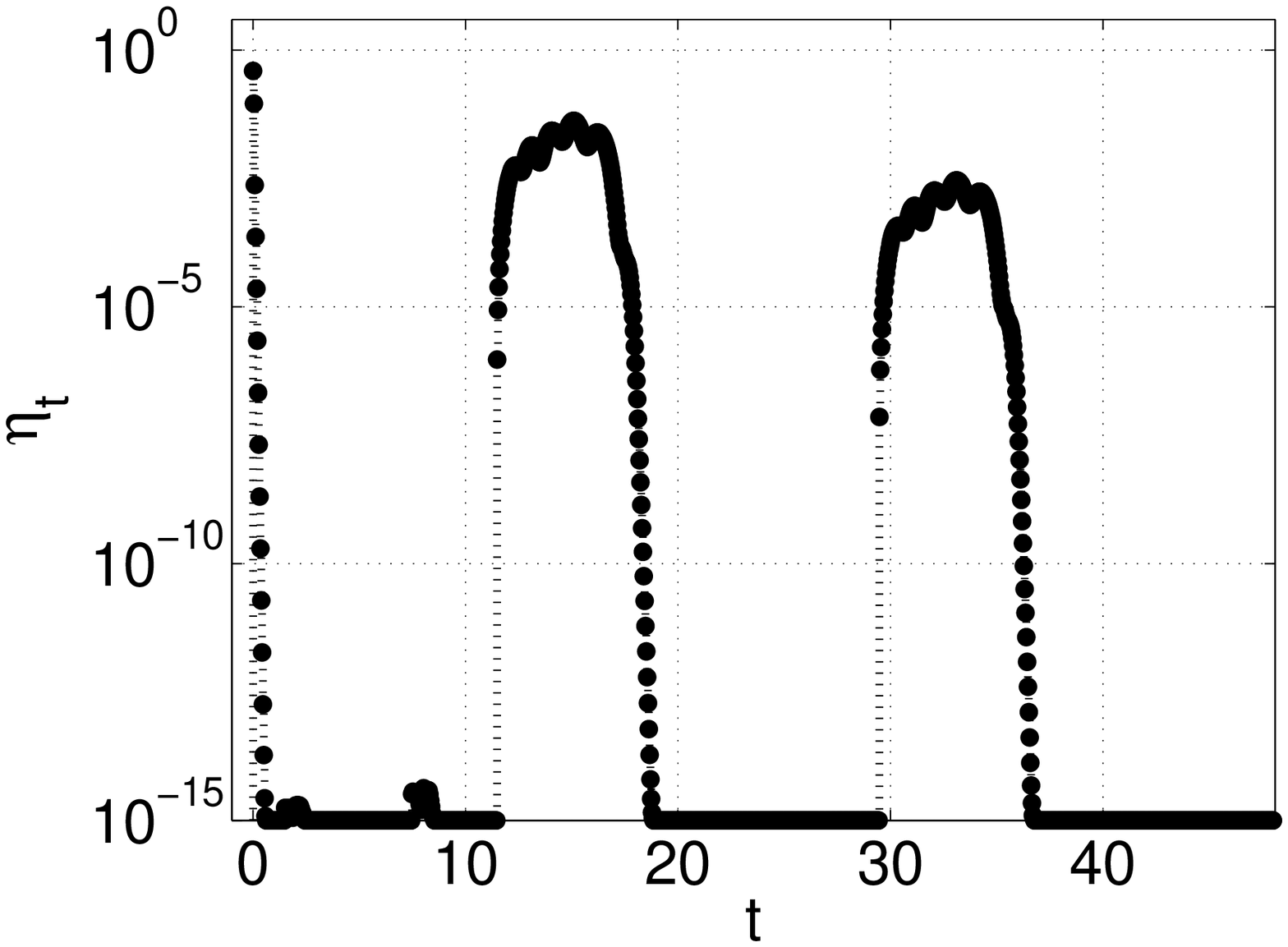, width=0.49\linewidth,height=0.16\textheight}}
}
\end{center}
\vspace{0.5cm}
\begin{center}{ 
 \subfigure[$\cfl(t_n)$, L=1, adaptive implicit timestep]
 	{\epsfig{file= 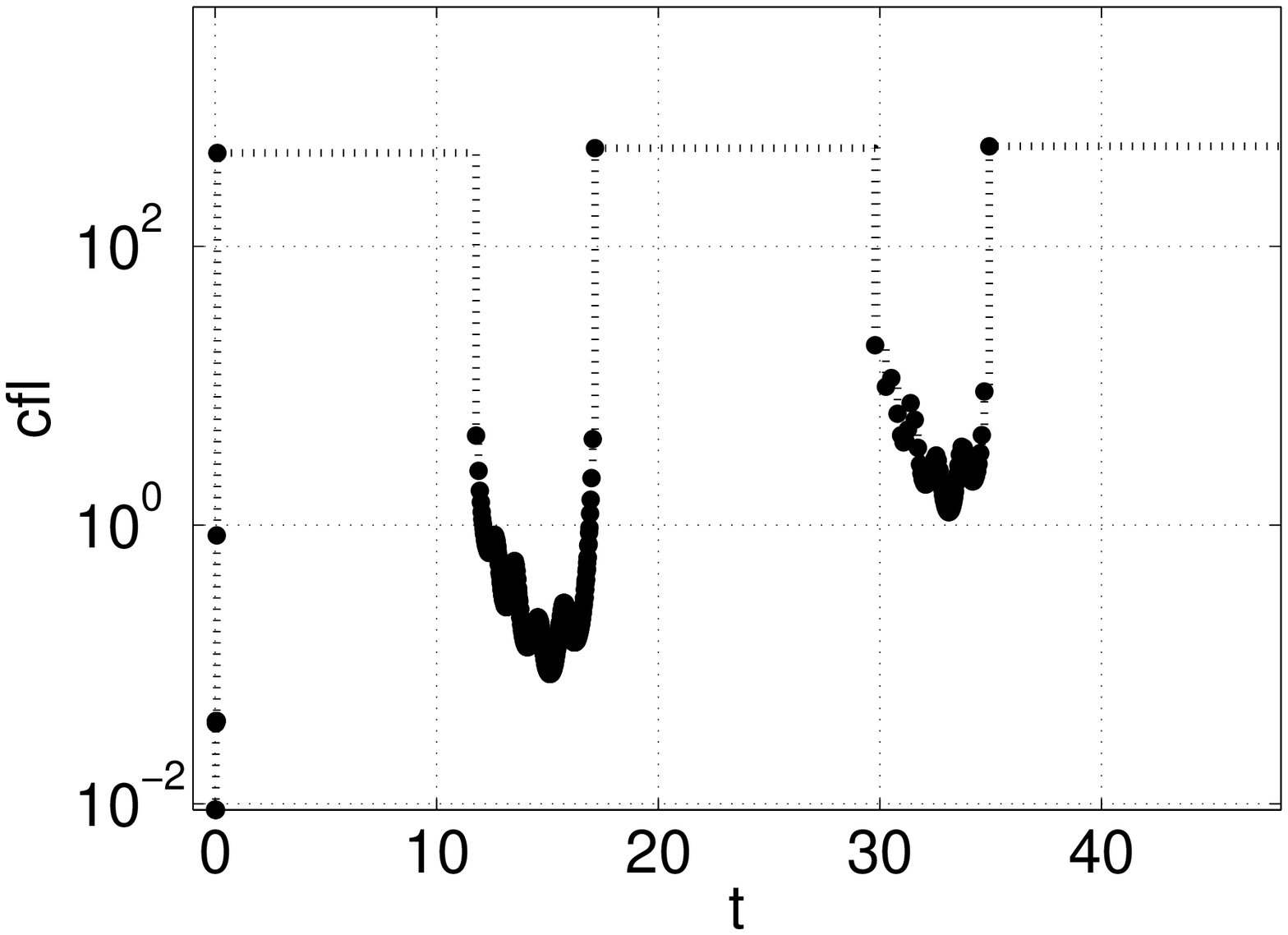, width=0.49\linewidth,height=0.16\textheight}}
 \subfigure[$\bar\eta_k^n(t_n)$, L=1, adaptive implicit timestep]
 	{\epsfig{file= 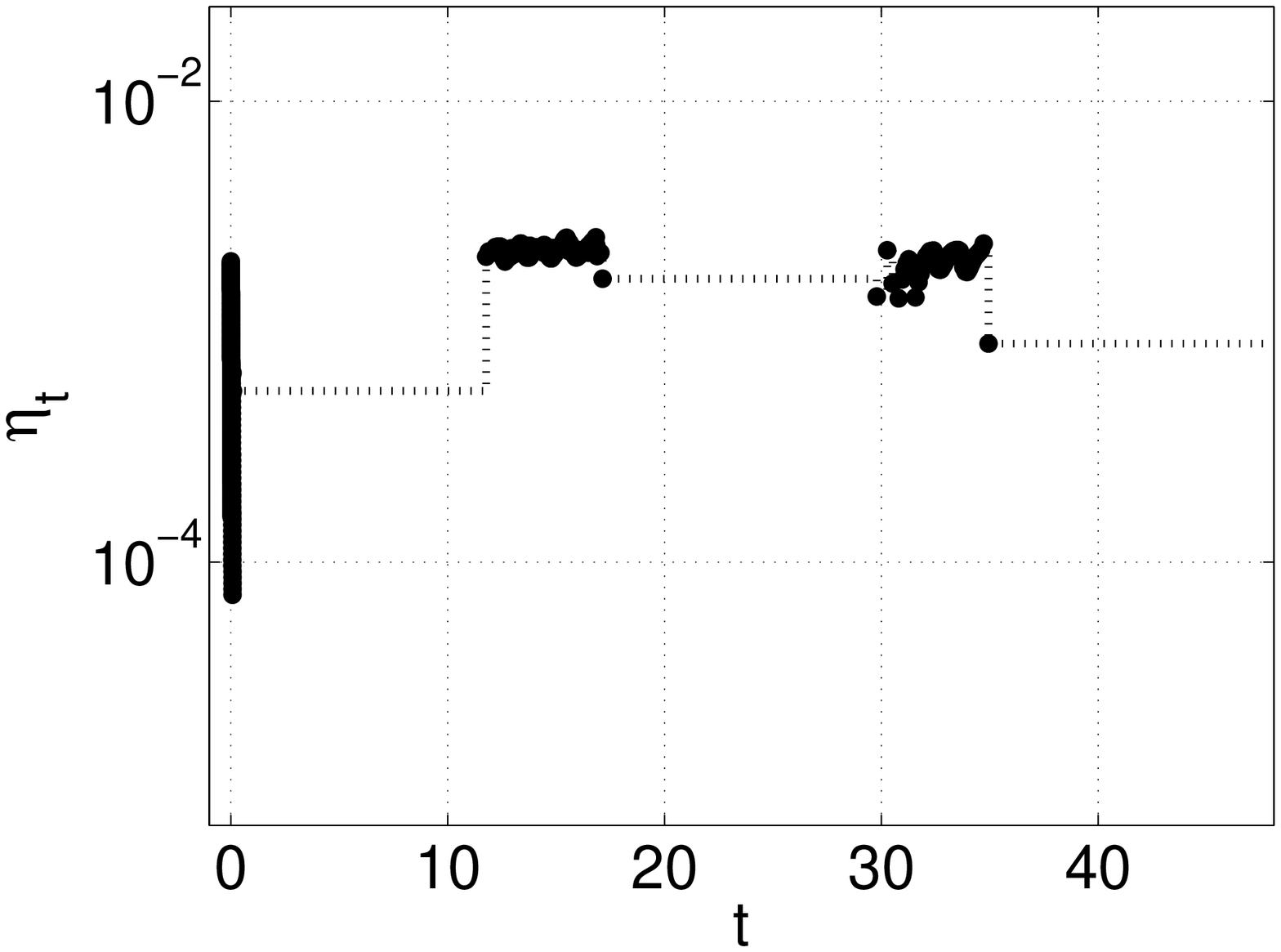,  width=0.49\linewidth,height=0.16\textheight}}
}
\end{center}
\caption{Example 1: Perturbed shock for Burgers' equation with equidistributed time error.
Left column: CFL(t); right column: $\bar\eta_k^n(t_n)$. Upper row: level $L=0$,
fully explicit scheme, uniform timestep. Lower row: level $L=1$,
fully implicit scheme, adaptive timestep.
(from top to bottom).} \label{fig.perturbed_shock.tol_0.01}
\end{figure}

Important additional information can be gained by looking at the plots in
Figure~\ref{fig.perturbed_shock.tol_0.01}, showing the CFL distribution on each
time interval $I_n$ and the normalized time components of the error estimator
$\ebkn$, both in logarithmic scale. The stationary and instationary regions are
separated by the estimator. In particular, note that
\begin{itemize}
\item
the time component of the error varies over more than 14 orders
of magnitude.
\item
in the three stationary regions, $\ebkn$ is very close to zero.
\item
the two instationary waves are distinguished very clearly. The second wave is about one order of magnitude smaller than the first wave. This corresponds closely to the different magnitudes of the inflow perturbations.
\item
furthermore, one can clearly identify an initial layer, where
$\ebkn=\mathcal O(1)$ at the inflow boundary $t=0$, and $\ebkn$
decays exponentially for time $t>0$ until it reaches machine accuracy.
\end{itemize}
We advance to level $L=1$, Figures~\ref{fig.perturbed_shock.tol_0.01}(c) and
(d). We observe that

\begin{itemize} 
\item
the error on level $L=1$ varies by less than 2 orders of magnitude, 12
orders of magnitude less than on level $L=0$. The magnitude of the
maximal error has decreased by almost two orders of magnitude.
Therefore the solution is much better resolved in the instationary regions,
and the computational recources are clearly distributed more efficiently.

\item
in the initial layer, the CFL number starts with $\mathcal{O}(10^{-2})$.
Then it grows at least exponentially until it reaches a maximal value of about
500. At the same time, the error $\bar\eta_k^n$ decays roughly by two orders of
magnitude. Thus, these initial steps can be seen as a preprocessing of the initial
data, to translate a prescribed steady shock on the pde level into a
steady discrete shock layer. The initial layer is also clearly visible
in the plot of the error distribution (and this will never disappear).
Indeed, the initial data, a sharp jump from 1 to -1, are a steady
shock only on the level of the exact solution. Numerically, the scheme has
to converge towards a discrete shock layer, and this will always need
a few time steps. In fact this is an instance of a scheme converging
towards a numerical steady state solution, using adaptive time steps.
\item
in the stationary region between the initial layer and the first
perturbation, the error is more than three orders of magnitudes
smaller than in the following flow field. Observe that the whole stationary region is
computed by a single time step. Thus the scheme is only held back from
choosing a larger time step by the appearance of the instationary
perturbation. If we had introduced this perturbation at a later time, the
time step and thus the local CFL number would haves been
correspondingly larger.
\item
The next region of stationary flow is again bridged by a single time step,
and correspondingly the local error is somewhat below the equidistributed
one.
\item
For the two perturbations, the normalized error is already close to being uniformly distributed.
\item 
Using large timesteps does not mean that each timestep has higher computational costs. Since the adaptation
 chooses large timesteps, where the solution is(nearly) stationary, these timesteps have low computational 
 costs.
\end{itemize} 
We would also like to point out one drawback of the equidistribution strategy
for the timestep. In the  first (and larger) instationary wave, the proposed CFL
number is often much smaller than unity, e.g.  $\min\limits_n(\cfl(t^n)) =
0.009$ in Figure~\ref{fig.perturbed_shock.tol_0.01}. It is well-known  that
lowering the CFL number much below unity smears the solution. Therefore, while 
such small timesteps may improve the temporal accuracy somewhat, they will deteriorate
the spatial accuracy considerably. Moreover, they increase the number of
timesteps, and hence the computational cost.  In the following example, we
discuss a more efficient strategy.

%%%%%%%%%% Example 2 %%%%%%%%%%%%%%%%
\vspace*{2ex}\noindent
{\bf Example 2:} 
This example is a modification of the first example which used a fully implicit
strategy for the timestep. Here we introduce a mixed {\em implicit/explicit strategy}. 
We still want the equidistribute the error, but we will give up this goal
partially when the local $\cfl$ number drops below a certain threshhold.

As discussed above, choosing timestep sizes with $\cfl$ much less than unity
seems to be inefficient both for  explicit and for implicit schemes. For
implicit methods, even timesteps with $\cfl <5$ are not efficient, since we have
to solve a nonlinear system of equations at each timestep. Thus, the new
implicit/explicit strategy switches to the cheaper (and less dissipative)
explicit method, if $\cfl < 5$, computing perhaps a few more timesteps if
$0.8<\cfl<5$, and saving timesteps if $\cfl<0.8$. (Of course, we could choose
other thresholds than $\cfl=0.8$ and 5.) 

As we can see in Table~\ref{table.perturbed_shock_cfl}, the new strategy
requires only 449 timesteps, instead of 1238 with the direct equidistribution
in Example 1. Out of these, only 166 are implicit and  hence expensive. This
leads to considerable speed-up. 

\begin{table}  
\begin{tabular}{c|c|c|c|c}
L	& N(expl) &$\bar\eta_k/J(u_h)$	&$\bar\eta_h/J(u_h)$	&$\bar\eta/J(u_h)$	\\ \hline
1	& 449\,(283) &8.50e-04		&2.14e-02	&2.23e-02
\end{tabular}\\  [1ex]
\caption{Example 2: Same as Table~\ref{table.perturbed_shock.1}, but CFL restriction from below
(implicit/explicit strategy).}
\label{table.perturbed_shock_cfl}
%\vspace*{-4ex}
\end{table}   
%
%
%timeit(45, 07759,0)
%tol =0.001
\begin{figure}[hbtp]
\begin{center}{ 
 \subfigure[$C\!F\!L(t_n)$, L=1, adapt.~impl./expl.~timestep]
 	{\epsfig{file= 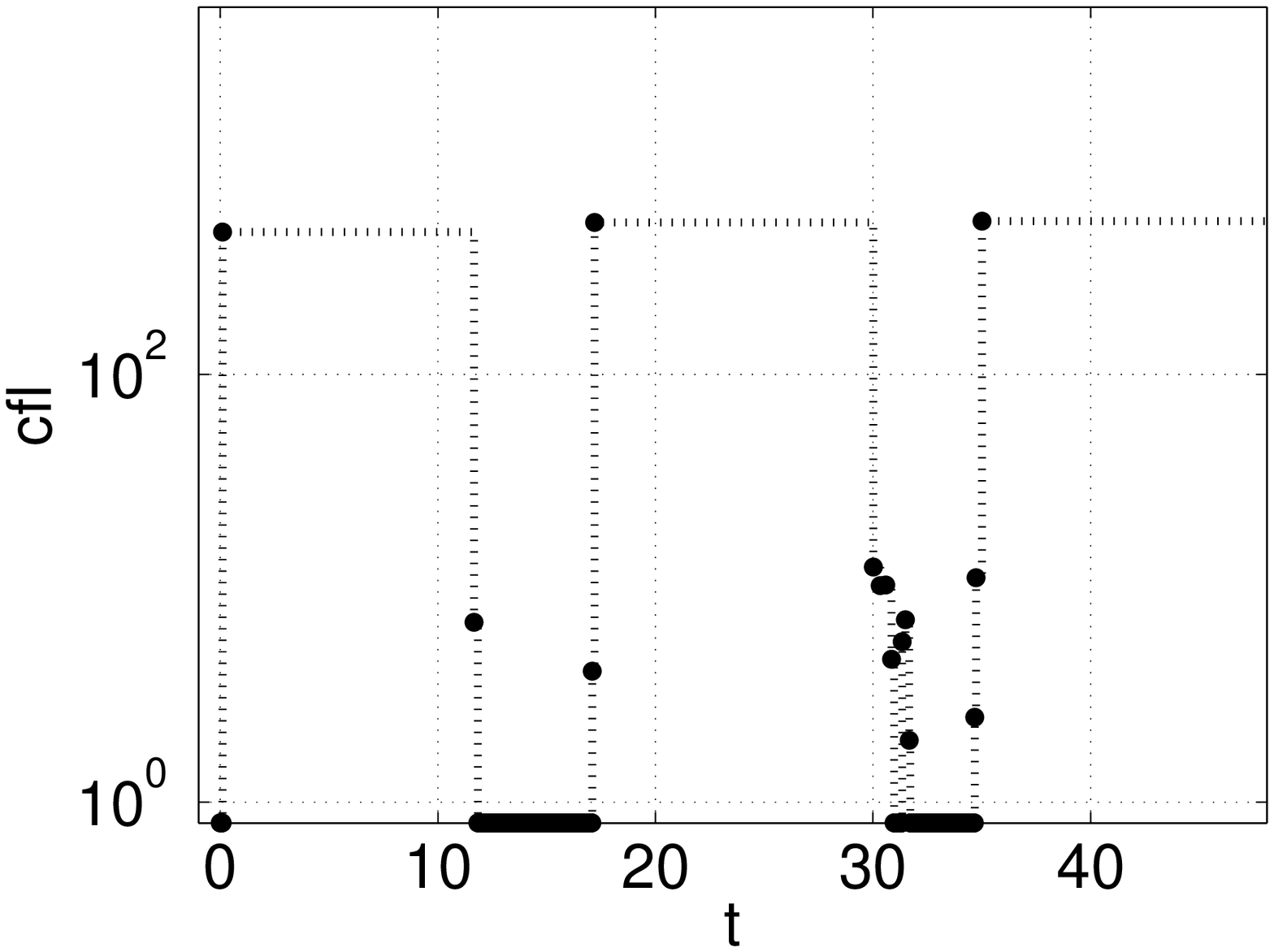, width=0.49\linewidth,height=0.16\textheight}}
 \subfigure[Error $\bar\eta_k^n$, L =1, adapt. impl./expl. timestep]
 	{\epsfig{file= 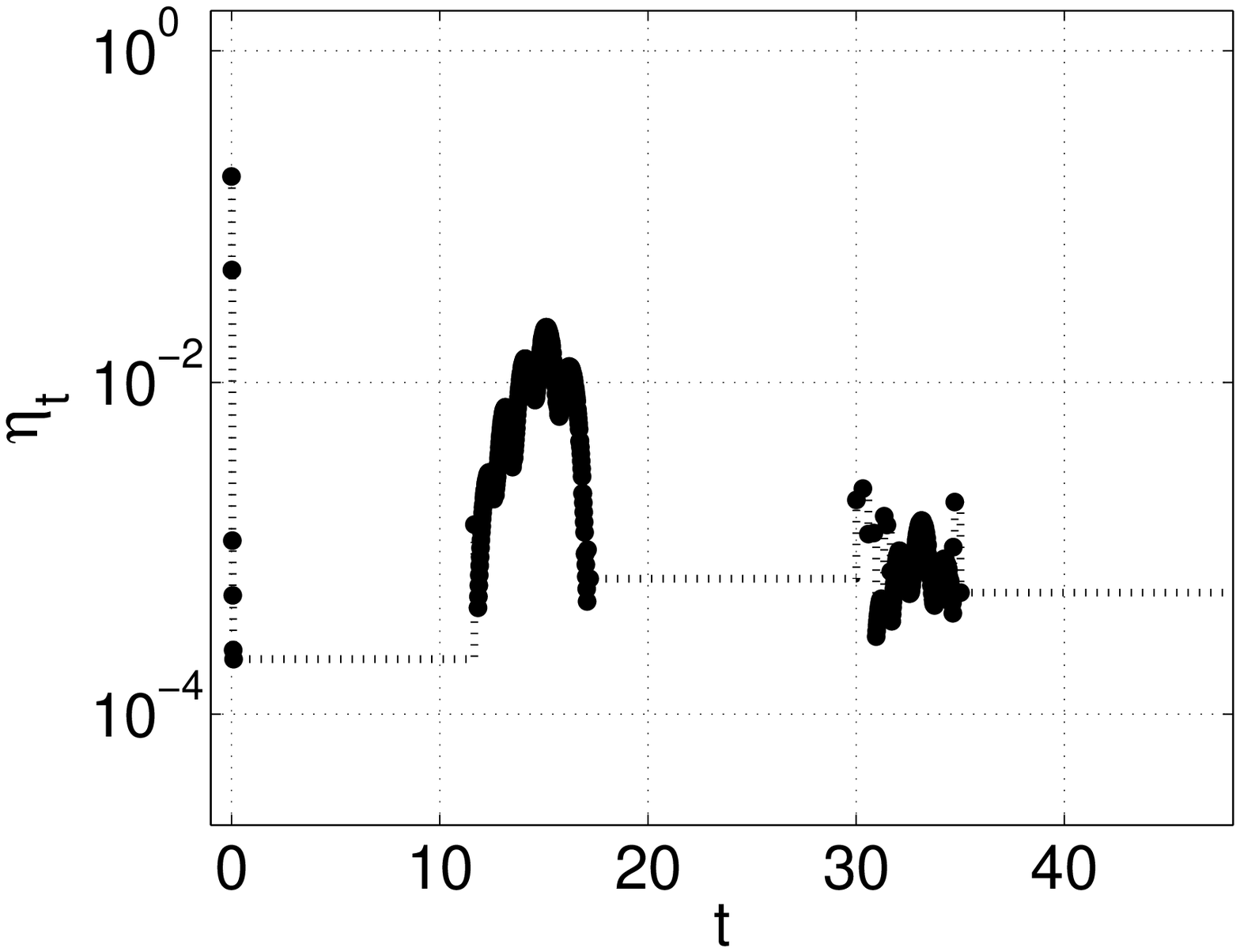,  width=0.49\linewidth,height=0.16\textheight}}
}
\end{center}
\caption{Example 2: Same as Figure~\ref{fig.perturbed_shock.tol_0.01} but
adaptive implicit/explicit strategy (CFL restriction from below).} 
\label{fig.perturbed_shock_cfl}
\end{figure} 

\begin{figure}[hbtp]
\begin{center} {
 \subfigure[$\cfl(t_n)$, uniform explicit]
 	{\epsfig{file= 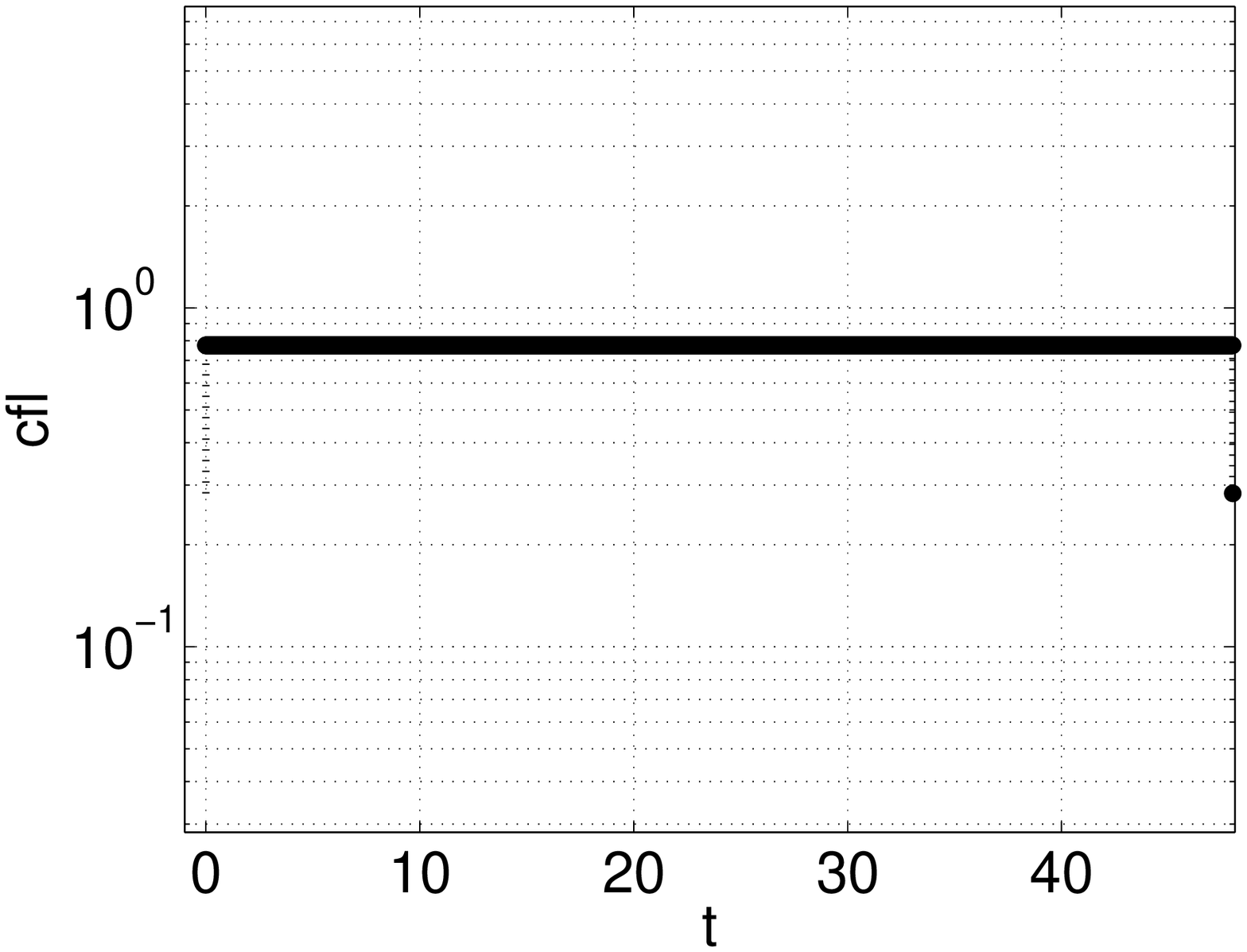, width=0.49\linewidth,height=0.16\textheight}}
 \subfigure[$\bar\eta_k^n(t_n)$, uniform explicit]
 	{\epsfig{file= 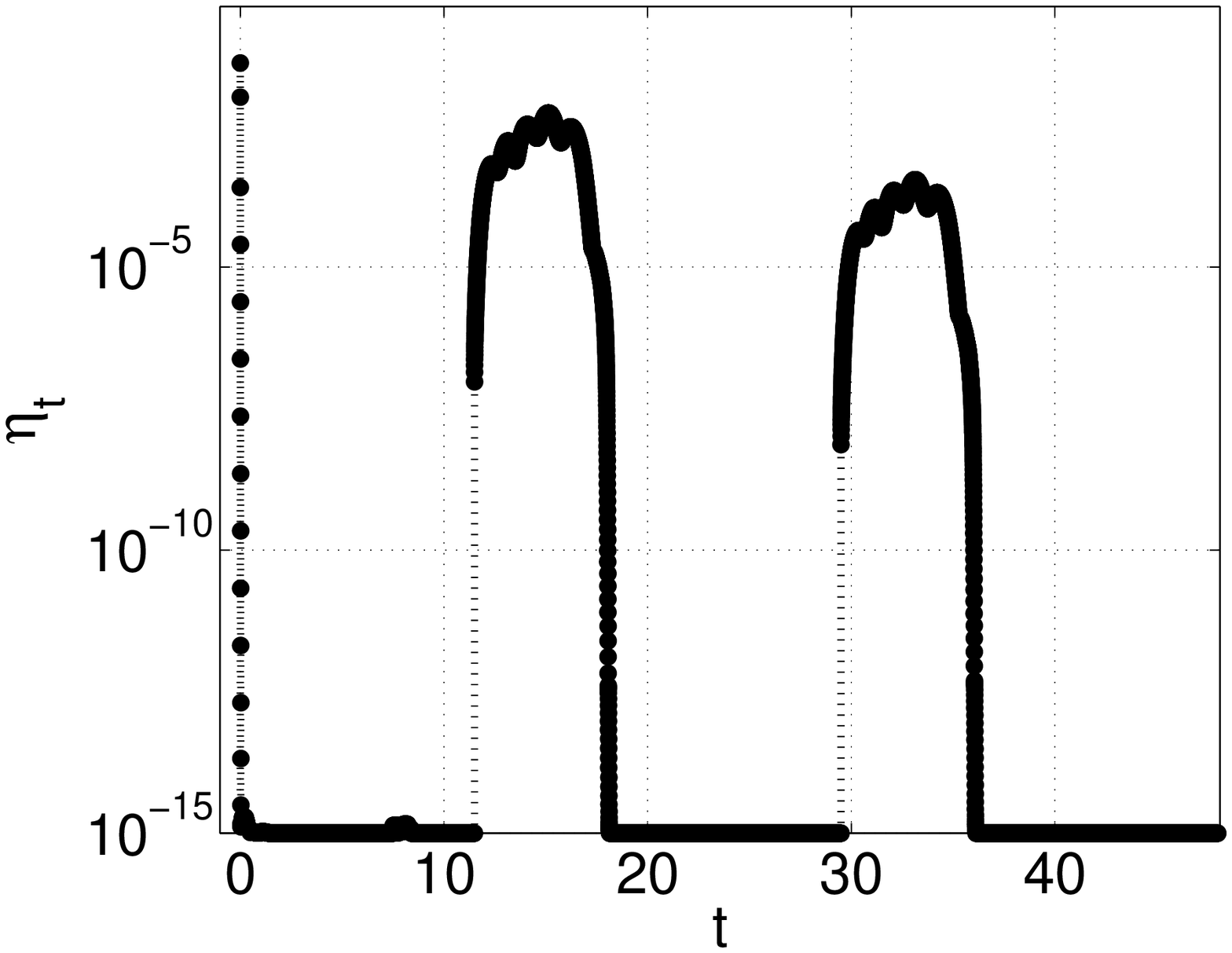,  width=0.49\linewidth,height=0.16\textheight}}
}
\end{center}
\vspace{0.5cm}
\begin{center}{ 
 \subfigure[$\cfl(t_n)$, adapt.~impl./expl., $\tol_k =2^{-4}\bar\eta_k^{ref}$]
 	{\epsfig{file= 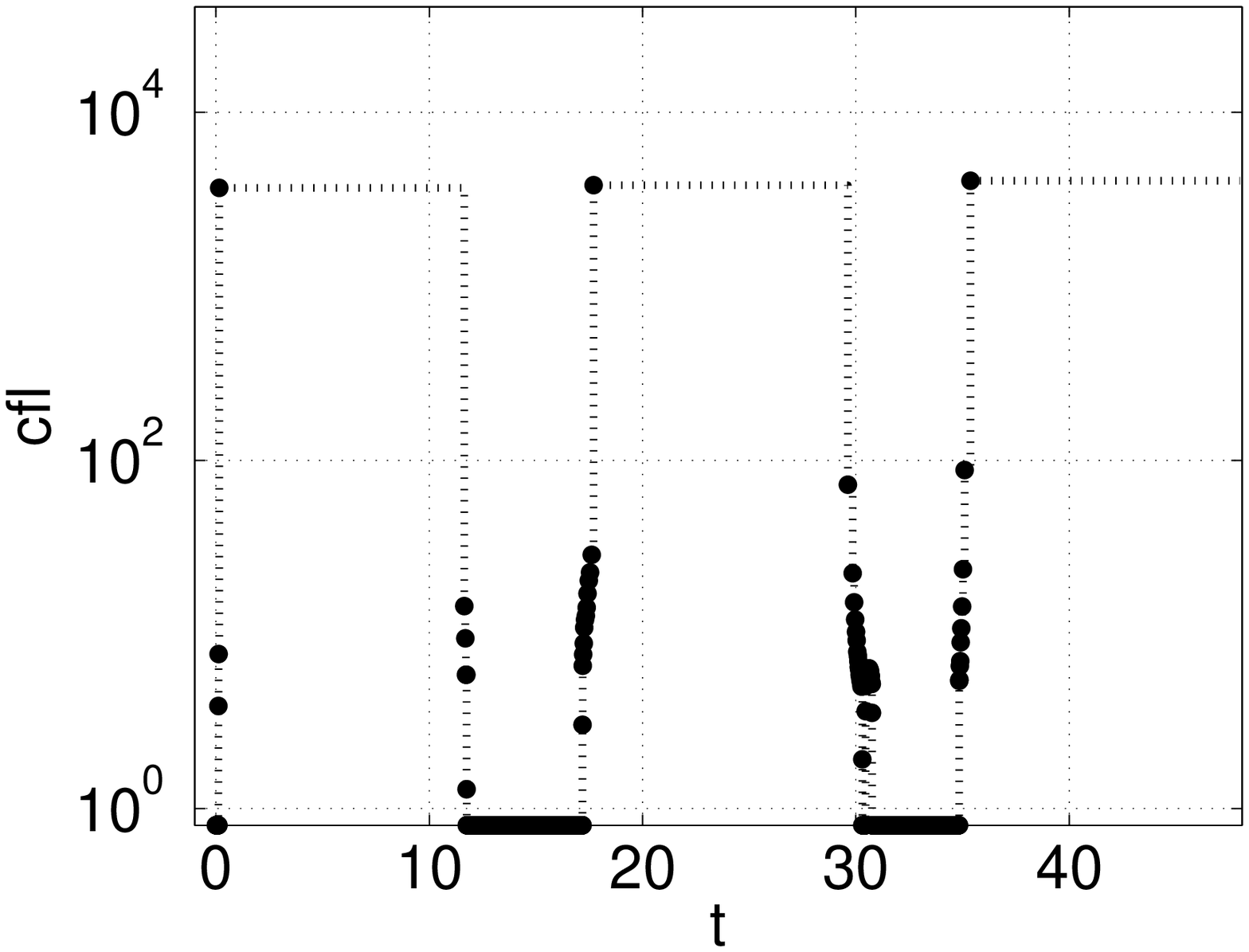, width=0.49\linewidth,height=0.16\textheight}}
 \subfigure[$\bar\eta_k^n(t_n)$, adapt.~impl./expl., $\tol_k =2^{-4}\bar\eta_k^{ref}$]
 	{\epsfig{file= 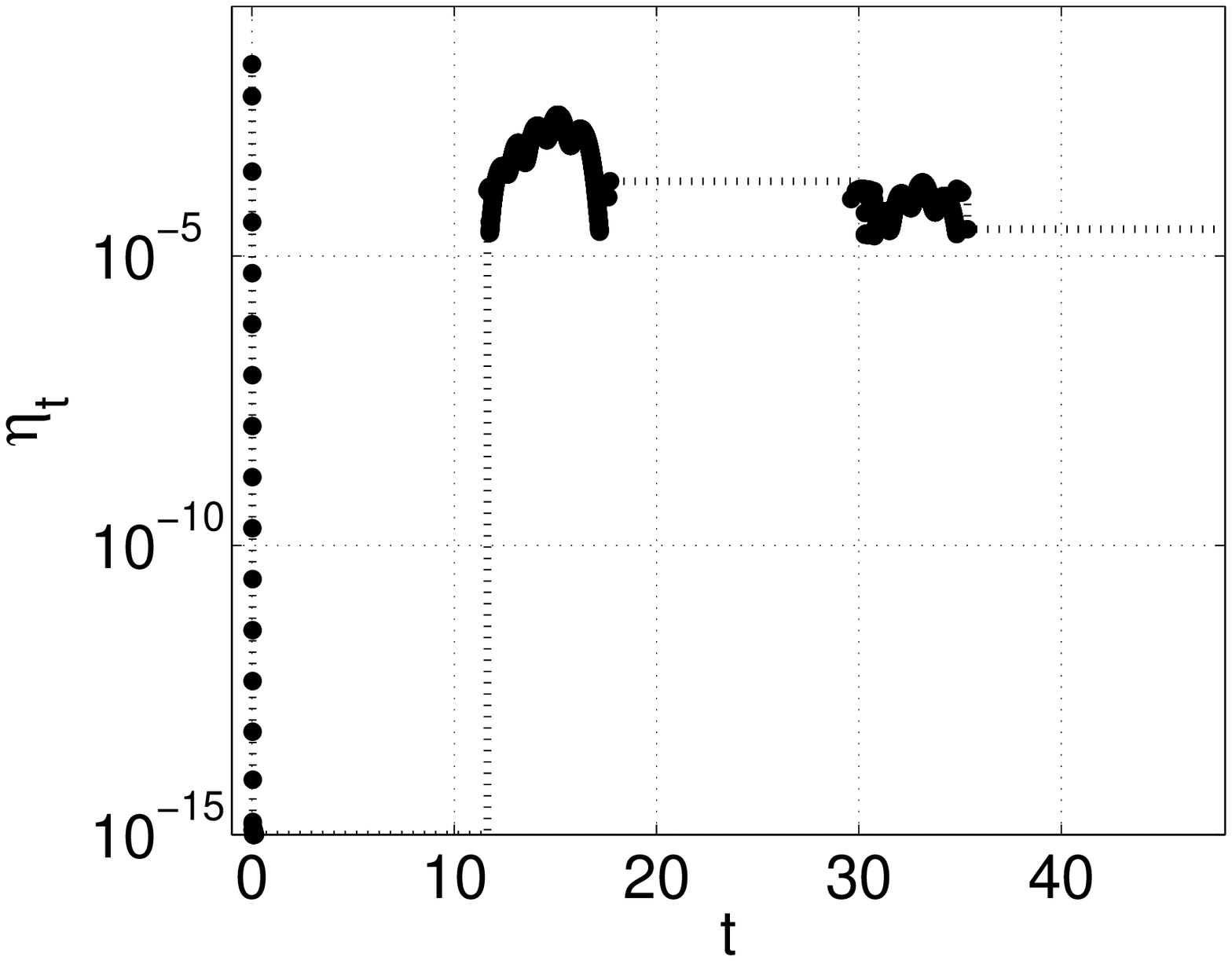,  width=0.49\linewidth,height=0.16\textheight}}
}
\end{center}
\vspace{0.5cm}
\begin{center}{ 
\subfigure[$\cfl(t_n)$, adapt.~impl./expl., $\tol_k =\bar\eta_k^{ref}$]
 	{\epsfig{file= 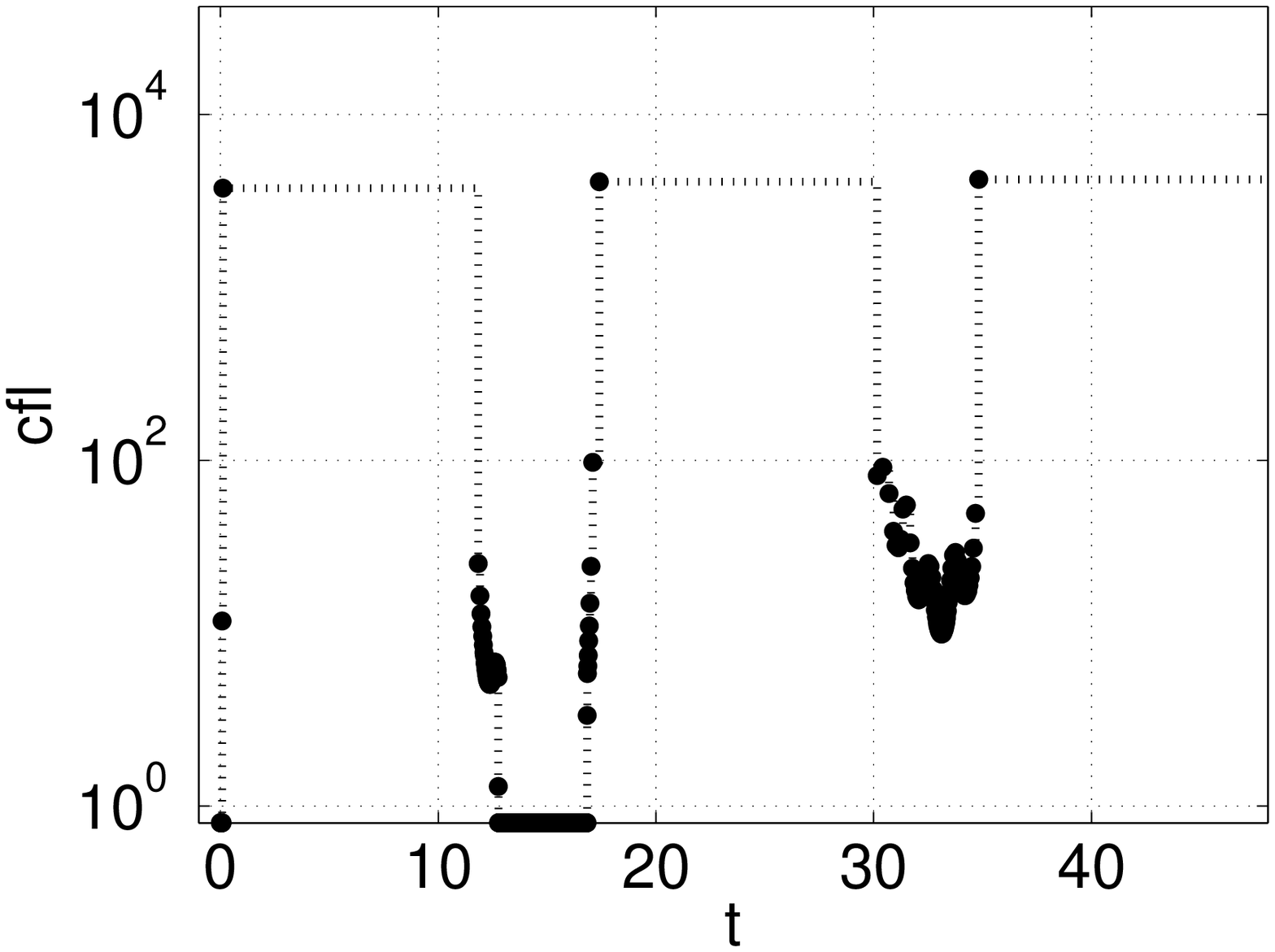, width=0.49\linewidth,height=0.16\textheight}}
 \subfigure[$\bar\eta_k^n(t_n)$, adapt.~impl./expl., $\tol_k =\bar\eta_k^{ref}$]
 	{\epsfig{file= 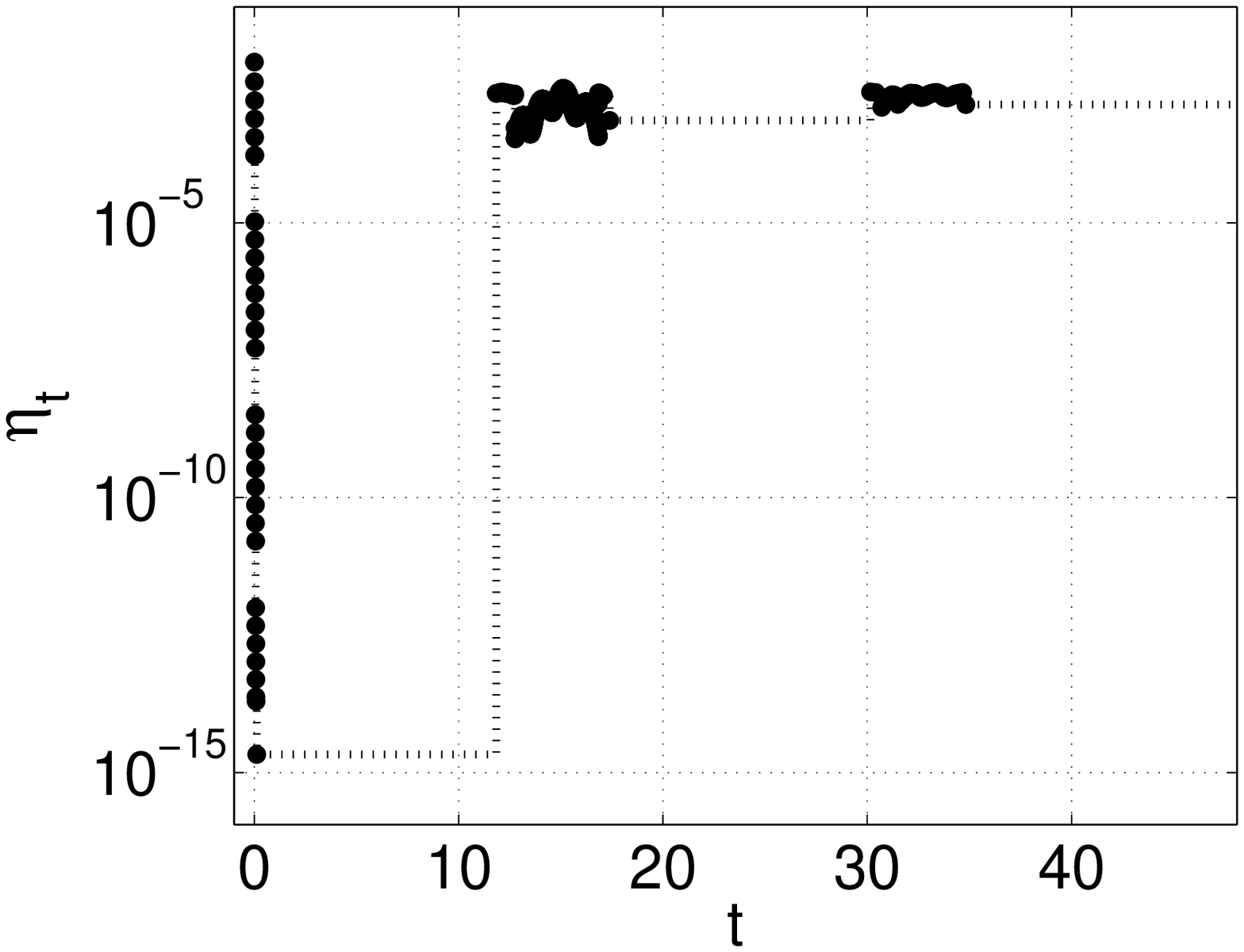, width=0.49\linewidth,height=0.16\textheight}}
}
\end{center}
\caption{Example 3: Same as Figure~\ref{fig.perturbed_shock_cfl} but on Level L =4 and with different tolerances. \label{fig.perturbed_shock_cfl_3}
}\end{figure}

%%%%%%%%%% Example 5 %%%%%%%%%%%%%%%% 
\vspace*{2ex}\noindent 
{\bf Example 3:}
Table~\ref{table.perturbed_shock_cfl_2} and Figure
~\ref{fig.perturbed_shock_cfl_3} show three extensions of Example 2. We used the
same implicit/explicit strategy as in Example 2, but after the explicit
reference computation on the coarse grid (L=0, error $\bar\eta_k^{ref}$), we
proceed directly to a finer grid with 320 cells (L=4). We compare an explicit
and two  implicit/explicit computations on the fine grid. 

The first row shows results of the fully explicit scheme with uniform refinement
in time and space for L=4. As expected, the errors are about $2^4$ times smaller
than those on the original coarse grid. Now suppose we wanted to reach
comparable errors on level $L=4$ using adaptive timestepping. Then we should set
the tolerance to be $\;\tol(4)=2^{-4}\bar\eta_k^{ref}\,$. The results of this
computation are shown in the second row of
Table~\ref{table.perturbed_shock_cfl_2}. The three components of the error are
comparable with those of the fully explicit computation, but the number of
timesteps is only 3975 instead of 19200. Out of these 3975 steps, only 1235 are
implicit.

Another strategy for equidistributing the error might be to fix any constant
tolerance, for example $\tol(4)=\bar\eta_k^{ref}$ itself. The results of this
computation are displayed in the last row of the table. The error in time is
now a factor 5-8 higher than for the other two computations, while the spatial
error is comparable. Remarkably, this computation needs only 1780 timesteps, and
only 507 of these are implicit. 

Both of these calculations show that considerable savings are possible with the
implicit/exp\-licit, time-adaptive strategy.

\begin{table}  
\begin{tabular}{c|c|c|c|c|c}
strategy	& $\tol_k$ 		& N(expl)			&$\bar\eta_k/J(u_h)$	&$\bar\eta_h/J(u_h)$	&$\bar\eta/J(u_h)$	\\ \hline
fully expl.	& --			& 19200\,(19200)		&7.11e-05	&4.57e-04		&5.28e-04	\\  
impl./expl.	& $2^{-4}\bar\eta_k^{ref}$	& \phantom{1}3975\,\phantom{1}(2740) 	&1.06e-04	&3.71e-04		&4.77e-04	\\ 
impl./expl.	& $\bar\eta_k^{ref}$		& \phantom{1}1780\,\phantom{1}(1243) 	&5.46e-04	&3.66e-04		&9.12e-04
\end{tabular}\\[1ex]  
\caption{Example 3: Same as Table~\ref{table.perturbed_shock_cfl}, but on level
L =4 and with different tolerances.}
\label{table.perturbed_shock_cfl_2}
\end{table}   

%
%

%%%%%%%%%%%%%%%%%%%%%%%%%%%%%%%%%%%%%%%%%%%%%%%%%%%%%%%%%%%%%%%%%%%%%%%%%%%%

\section{Conclusions}
\label{sec.conclusions}

In this paper, we combine space- and time-projections of S\"uli, Houston and Hartmann
to split the classical adjoint based error representation formula for
target functionals into space and time components. Based on a numerical study
of these components we design an adaptive strategy which attempts to minimize the number of
time steps by equidistributing the time components of the error. We apply the adaptive
scheme to a weak perturbation of a stationary shock.

Already on a very coarse mesh of 20 points the error representation formula
precisely gives the location and strength of the instationary perturbations.
This can be translated into efficient timestep distributions, which respect a
desired accuracy. We show that these timestep distributions can be applied
successfully to much finer spatial grids. 

We never compute implicit timesteps below CFL=5. Instead, when the error
analysis suggests a  timestep below CFL=5, we switch to an explicit scheme with
CFL=0.8. This implicit/explicit strategy gives considerable savings.

For nonlinear perturbations of a stationary shock, we have demonstrated that our
strategy does reach its goals: it separates initial layers, stationary regions
and perturbations cleanly and chooses just the right timestep for each of them.

Besides building upon well-established adjoint techniques, we have also added
a new ingredient which simplifies the computation of the dual problem. We show
that it is sufficient to compute the  spatial gradient of the dual solution,
$w=\nabla \varphi$, instead of the dual solution $\varphi$ itself. This gradient
satisfies a conservation law instead of a transport equation, and it can
therefore be computed with the same algorithm as the forward problem, and in the
same finite element space. For discontinuous transport coefficients, the new
conservative algorithm for $w$ is more robust than our previous transport
schemes for $\varphi$. 

In ongoing work, we are adapting this strategy to aerodynamic 
problems. First test calculations show a promising speed-up.

%%%%%%%%%%%%%%%%%%%%%%%%%%%%%%%%%%%%%%%%%%%%%%%%%%%%%%%%%%%%%%%%%%%%%%%%%%%%

\medskip

{\bf Acknowledgement:} The authors would like to thank Ralf Hartmann, Paul
Houston, Mario Ohlberger and Endre S\"uli for stimulating discussions.
The work of both authors was supported by DFG grant 
SFB 401 at RWTH Aachen. Part of the work was completed while the first author
was in residence at the Center of Mathemtics for Applications (CMA) at Oslo University. Both authors would like to thank
CMA its members for their generous hospitality.
 
%%%%%%%%%%%%%%%%%%%%%%%%%%%%%%%%%%%%%%%%%%%%%%%%%%%%%%%%%%%%%%%%%%%%%%%%%%%%%

%%%%%%%%%%%%%%%%%%%%%%%%%%%%%%%%%%%%%%%%%%%%%%%%%%%%%%%%%%%%%%%%%%%%%%%%%%%%%%%

%\include{answer}
%%%%%%%%%%%%%%%%%%%%%%%%%%%%%%%%%%%%%%%%%%%%%%%%%%%%%%%%%%%%%%%%%%%%%%%%%%%%%%%
\end{document}